\documentclass[12pt,a4paper]{article}
\usepackage{amssymb,amsmath,amsthm,graphicx,verbatim,url,color} 
\usepackage{cancel}
\usepackage{hyperref}
\newtheorem{theorem}{\sc Theorem.}[section]
\newtheorem{lemma}[theorem]{\sc Lemma.}
\newtheorem{remark}[theorem]{\sc Remark.}

\renewcommand{\theequation}{\arabic{section}.\arabic{equation}}

\newcommand{\bRgeq}{{\mathbb R}_{\geq 0}}

\newcommand{\bR}{{\mathbb R}}
\newcommand{\bS}{{\mathbb S}}

\newcommand{\nabs}{\nabla_{\!s}}
\newcommand{\matpartu}{\partial_t^\bullet}

\newcommand{\matpartx}{\partial_t^\circ}
\newcommand{\matpartxh}{\partial_t^{\circ,h}}
\newcommand{\Id}{I\!d}
\newcommand{\bigchi}{\ensuremath{\mathrm{\mathcal{X}}}}
\newcommand{\charfcn}[1]{\bigchi_{#1}} 
\newcommand{\subT}{{\rm T}}            
\newcommand{\subN}{{\rm N}}
\newcommand{\subH}{{\rm B}}
\newcommand{\Gauss}{{\mathcal{K}}}
\newcommand{\onbtau}{\vec{\mathfrak t}}
\def\epsilon{\varepsilon} 
\def\hat{\widehat}

\newcommand{\dH}[1]{\;{\rm d}{\mathcal{H}}^{#1}} 
\newcommand{\Wht}{S^h_{\ell-1}(\Gamma^h(t))}

\newcommand{\qVhz}{[S^h_\ell(\Gamma^0)]^3}   
\newcommand{\qVh}{[S^h_\ell(\Gamma^m)]^3}
\newcommand{\qVht}{[S^h_\ell(\Gamma^h(t))]^3}

\newcommand{\Wh}{S^h_{\ell-1}(\Gamma^m)}
\newcommand{\uspace}{\mathbb{U}}

\newcommand{\GT}{{\mathcal{G}_T}}
\newcommand{\GhT}{{\mathcal{G}^h_T}}
\newcommand{\ddt}{\frac{\rm d}{{\rm d}t}}
\newcommand{\id}{{\rm id}}

\newcommand{\mat}[1]{\underline{\underline{#1}}\rule{0pt}{0pt}}
\begin{document}
\title{
A parametric finite element method for the incompressible
Navier--Stokes equations on \\ an evolving surface}
\author{Harald Garcke\footnotemark[2]\ \and 
        Robert N\"urnberg\footnotemark[3]}

\renewcommand{\thefootnote}{\fnsymbol{footnote}}
\footnotetext[2]{Fakult{\"a}t f{\"u}r Mathematik, Universit{\"a}t Regensburg, 
93040 Regensburg, Germany \\ {\tt harald.garcke@ur.de}}
\footnotetext[3]{Dipartimento di Mathematica, Universit\`a di Trento,
38123 Trento, Italy \\ {\tt robert.nurnberg@unitn.it}}

\date{}

\maketitle

\begin{abstract}
In this paper we consider the numerical approximation of the
incompressible surface Navier--Stokes equations on an evolving surface.
For the discrete representation of the moving surface we use 
parametric finite elements of degree $\ell \geq 2$. In the semidiscrete
continuous-in-time setting we are able to prove a stability estimate that
mimics a corresponding result for the continuous problem. Some numerical
results, including a convergence experiment, demonstrate the practicality and
accuracy of the proposed method.
\end{abstract} 
\renewcommand{\thefootnote}{\arabic{footnote}}

\section{Introduction} \label{sec:1}

We consider an evolving closed, compact and oriented surface $\Gamma$ in $\bR^3$
whose evolution is driven by the velocity $\vec u$
that satisfies the incompressible surface Navier--Stokes equations
\begin{equation} \label{eq:esns}
\rho \, \matpartu\,\vec u - \nabs\cdot\mat\sigma =
\vec g + \alpha\, f_\Gamma\,\vec \nu,\qquad
\nabs\cdot\vec u = 0,\qquad \mathcal V = \vec u\cdot\vec\nu.
\end{equation}
Here $\matpartu$ denotes the material time derivative on $\Gamma$,
$\nabs\cdot$ is the surface divergence,
\begin{equation} \label{eq:BStensor}
\mat\sigma= 2\,\mu\,\mat D_s (\vec u)  - p\,\mat{\mathcal{P}}
\end{equation}
is the surface stress tensor, $\vec\nu$ is a unit normal on $\Gamma$, 
$\mathcal{V}$ is the normal velocity of the evolving surface $\Gamma$,
$\vec g$ is some external forcing
and $\rho$ is the constant density. In addition,
 $\mu \in \bRgeq$ is the interfacial shear viscosity, $D_s (\vec u)$ is the surface rate-of-deformation tensor, $\mat{\mathcal{P}}$ is the projection on the tangent space and
$p$ denotes the surface pressure, which acts as a Lagrange
multiplier for the incompressibility condition. 
For the forcing $f_\Gamma$ we have some
intrinsic properties of the surface in mind, and $\alpha \in \bRgeq$ is
a related coefficient. The second condition in \eqref{eq:esns} models the
incompressibility of the surfaces material, meaning in particular that the
total surface area is conserved. The third equation in \eqref{eq:esns} means
that the surface evolves according to the normal component of the velocity
$\vec u$.

For example, the forcing $f_\Gamma$ in \eqref{eq:esns} can be derived from
a simple bending energy, e.g.\
\begin{equation} \label{eq:bendingE}
E(\Gamma) = \tfrac12 \int_{\Gamma} \varkappa^2 \dH{2}\,.
\end{equation}
By $\varkappa=-\nabs\cdot\vec \nu$   we denote the mean curvature 
(the sum of the principal curvatures) of $\Gamma$,
with the sign convention for the normal that $\varkappa=-2$ for the unit sphere,
and $\dH{2}$ 
indicates integration with respect to the two-dimensional surface measure. 
The forces from the bending energy act in a direction normal to the surface 
and $f_\Gamma $ is given as minus the first variation of $E(\Gamma)$, i.e.,
\begin{equation}\label{eq:fGamma}
f_\Gamma = -\Delta_s\,\varkappa
-\varkappa \,|\nabs\,\vec \nu|^2
+\tfrac{1}{2}\,\varkappa^3\,, 
\end{equation}
where $\Delta_s$ is the surface Laplacian and $\nabs$ is the surface
gradient. For a derivation of the $f_\Gamma$-term, we refer to Sections~9 and 10 in \cite{bgnreview}. 

Interfacial fluid mechanics was first thoroughly discussed by
\cite{Scriven60}, generalizing 
earlier ideas of Boussinesq. In this context the surface
stress tensor \eqref{eq:BStensor} was first introduced, and is hence
often called the Boussinesq--Scriven tensor.
The model \eqref{eq:esns} has been derived first in curvilinear coordinates in \cite{HuZE07} using 
conservation laws of surface mass and momentum quantities 
as physical principles. The same conservation laws have been used in \cite{JankuhnOR18} 
to derive surface Navier–Stokes equations in Cartesian coordinates in
$\bR^3$. It is also possible to derive the surface Navier--Stokes equations in $\bR^3$ using, instead of a
surface momentum conservation principle, a variational energy principle, which has been done in \cite{KobaLG17}.
In \cite{Miura18} the physical conservation laws of volume mass and momentum quantities are used
 in a thin tubular neighbourhood of an evolving surface. Using a thin film limit the author then derives the 
surface Navier–Stokes equations on an evolving hypersurface in $\bR^3$.
Similar as in the thin film derivation just mentioned, \cite{NitschkeRV20} derived tangential surface Navier–Stokes equations 
formulated with the help of a curvilinear coordinate system.

The model \eqref{eq:esns} can be viewed as a simplified variant of a model for fluidic membranes studied in \cite{ArroyoS09,nsns} which couples the bulk Navier--Stokes equations
to surface Navier--Stokes equations and also involves a bending energy.
The model \eqref{eq:esns} with $\vec u \cdot \vec \nu = V$ a given 
velocity is considered in \cite{OlshanskiiRZ22,OlshanskiiRS24}.
Frequently also Navier--Stokes equations on stationary surfaces have been considered and we refer to 
\cite{Fries18,OlshanskiiRZ21,PrussSW21,SimonettW22} 
for papers dealing with this setting.
For a nice review on Navier--Stokes equations on evolving surfaces, we refer to \cite{BrandnerRS22}.
Well-posedness results for \eqref{eq:esns} with \eqref{eq:fGamma} have been shown in \cite{WangZZ12,AbelsLP25preprint}.

Finally, let us mention that
for the energy \eqref{eq:bendingE}, and for $\rho=0$, the resulting evolution 
law \eqref{eq:esns} can be considered as a purely geometric evolution of the 
surface, which is close to Willmore flow with surface area conservation. 
If volume conservation is also considered, then the related flow is Helfrich 
flow, see e.g. \cite{willmore,pwf,pwfade,BonitoNP10,NeunteufelSS23,BrazdaKS24}.

Of interest is also the case $f_\Gamma = 1$ and then $\alpha$ can play the role of a Lagrange
multiplier for the side constraint 
$\int_\Gamma \mathcal{V}\dH{2} = \int_\Gamma \vec u\cdot\vec\nu \dH{2} =0$,
which would enforce the conservation of the volume enclosed by $\Gamma$.
We also observe that replacing \eqref{eq:bendingE} with the total surface area
$\int_\Gamma 1 \dH{2}$ would result in the
variation $f_\Gamma=\varkappa$. However, using that forcing in \eqref{eq:esns}
would have no effect on the flow, due to the surface's incompressibility. In
particular, the solution to \eqref{eq:esns} would be independent of the value
of $\alpha$. 

There has been considerable interest in the numerical approximation of solutions to \eqref{eq:esns} and its variants. These variants include, for example, evolving manifolds where the normal velocity is prescribed, as well as stationary surfaces. Some use parametric finite elements 
\cite{NitschkeVW12,ReutherV15,NitschkeRV17,ReutherV18,Fries18,ReutherNV20,BrandnerJPRV22,KrauseKV23,KrauseV23,ElliottS25}
and others use TraceFEM involving ideas from level set methods, see e.g.\ 
\cite{JankuhnORZ21,OlshanskiiRZ22,OlshanskiiRS24}.
However, to the best knowledge of the authors, there appears to be no numerical method for the full Navier--Stokes equations on an evolving surface, i.e.\ for the system \eqref{eq:esns} with either $\alpha=0$
or with $\alpha>0$, for which one can prove stability estimates.
It is the main goal of this paper to for the first time introduce such a method, to rigorously show stability estimates for a semidiscrete variant and to show via numerical experiments that the method works well in practice.

The outline of the paper is as follows. In Section~\ref{sec:gov} we introduce the governing equations. In Section~\ref{sec:weak}, we derive a weak formulation which   is the basis for a semidiscrete finite element discretization which is derived in Section \ref{sec:sd}. For the resulting set of equations, we can show stability estimates.
The fully  discrete finite element method and the solution techniques for the resulting linear systems are discussed in 
Sections~\ref{sec:fd} and \ref{sec:sol}. We also show that the discrete systems possess a unique solution provided that a
Ladyzhenskaya--Babu\v{s}ka--Brezzi (LBB) condition holds.
In Section~\ref{sec:nr2} we present results from numerical computations 
including tables reporting on errors and the order of convergence as well as computations showing qualitative properties of the solution.
We end the paper with an appendix in which we prove formulas which are important to separate the surface Navier--Stokes system in a normal and a tangential part, and in which we derive also an explicit radially symmetric solution.

\setcounter{equation}{0}
\section{Governing equations} \label{sec:gov}
We assume that $(\Gamma(t))_{t\in [0,T]}$ 
is a sufficiently smooth evolving hypersurface without boundary that is
parameterized by $\vec x(\cdot,t):\Upsilon\to\bR^3$,
where $\Upsilon\subset \bR^3$ is a given reference manifold, i.e.\
$\Gamma(t) = \vec x(\Upsilon,t)$. Then
\begin{equation*} 
\vec{\mathcal{V}}(\vec z, t) = \vec x_t(\vec q, t)
\qquad \forall\ \vec z = \vec x(\vec q,t) \in \Gamma(t)
\end{equation*}
defines the velocity of $\Gamma(t)$, and
$\mathcal{V} = \vec{\mathcal{V}} \cdot\vec{\nu}$ is
the normal velocity of the evolving hypersurface $\Gamma(t)$,
where $\vec\nu(t)$ is a unit normal on $\Gamma(t)$.
Moreover, we define the space-time surface
$\GT = \bigcup_{t \in [0,T]} \Gamma(t) \times \{t\}$.
Throughout this paper, for notational convenience, 
we often identify $\Gamma(t) \times \{t\}$ with $\Gamma(t)$.

On the free surface $\Gamma(t)$, the following conditions need to hold:
\begin{subequations} \label{eq:1}
\begin{alignat}{2}
\rho \, \matpartu\,\vec u - \nabs\cdot\mat\sigma & =
\vec g + \alpha\, \vec f_\Gamma 
\qquad &&\mbox{on } \Gamma(t)\,, \label{eq:1b} \\ 
\nabs\cdot\vec u & = 0
\qquad &&\mbox{on } \Gamma(t)\,, \label{eq:1c} \\ 
\vec{\mathcal{V}}\cdot\vec\nu &= 
\vec u\cdot\vec \nu \qquad &&\mbox{on } \Gamma(t)\,, 
\label{eq:1d} 
\end{alignat}
\end{subequations}
where $\rho \in \bRgeq$
denotes the surface material density,
$\alpha \in \bRgeq$ is a coefficient and 
$\vec f_\Gamma = f_\Gamma\,\vec\nu$ is defined by e.g.\ (\ref{eq:fGamma}).
In addition, 
$\nabs\cdot$ denotes the surface divergence on $\Gamma(t)$, and the surface
stress tensor is given by
\begin{equation} \label{eq:sigmaG}
\mat\sigma = 2\,\mu\,\mat D_s (\vec u) 
 - p\,\mat{\mathcal{P}}
\quad\text{on}\ \Gamma(t)\,,
\end{equation}
where $\mu \in \bRgeq$ is the interfacial shear viscosity and
$p$ denotes the surface pressure, which acts as a Lagrange
multiplier for the incompressibility condition (\ref{eq:1c}). 
Here
\begin{subequations}
\begin{equation} \label{eq:Ps}
\mat{\mathcal{P}} = \mat\Id - \vec \nu \otimes \vec \nu
\quad\text{on}\ \Gamma(t)
\end{equation}
is the projection onto the tangent space of $\Gamma(t)$, and
\begin{equation} \label{eq:Ds} 
\mat D_s(\vec u) = \tfrac12\,\mat{\mathcal{P}}\,(\nabs\,\vec u + 
(\nabs\,\vec u)^T)\, \mat{\mathcal{P}}
\quad\text{on}\ \Gamma(t)\,,
\end{equation}
\end{subequations}
is the surface rate-of-deformation tensor.
Moreover, $\nabs = \mat{\mathcal{P}} \,\nabla = 
(\partial_{s_1}, \partial_{s_2}, \partial_{s_3})$ 
denotes the surface gradient on $\Gamma(t)$, and
$\nabs\,\vec u = \left( \partial_{s_j}\, u_i \right)_{i,j=1}^3$.
We recall from Lemmas~13 and 14 in \cite{bgnreview} 
that on writing 
$\vec\xi= \xi_N\,\vec\nu + \vec\xi_\subT$, 
with $\vec\xi_\subT = \mat{\mathcal{P}}\,\vec\xi$, it holds that 
\begin{subequations} \label{eq:diff2}
\begin{equation} \label{eq:Ds2}
\mat D_s(\vec\xi) = \xi_\subN\, \nabs\,\vec\nu
+ \tfrac12\,(\mat{\mathcal{P}}\,\nabs\,\vec\xi_\subT + 
(\nabs\,\vec\xi_\subT)^T\mat{\mathcal{P}}),
\end{equation}
as well as
\begin{equation} \label{eq:nabs2}
\nabs \cdot \vec\xi = \xi_\subN\, \nabs\cdot\vec\nu
+ \nabs\cdot \vec\xi_\subT = - \xi_\subN\, \varkappa
+ \nabs\cdot \vec\xi_\subT.
\end{equation}
\end{subequations}
It is clear from 
\eqref{eq:diff2} that on parts of $\Gamma(t)$ with vanishing principal curvatures, i.e., 
where the surface is flat, controlling the norms of the
differential operators on the left hand sides offers no control on
$\xi_\subN=\vec\xi\cdot\vec\nu$. This has a profound impact on the analysis of
(Navier--)Stokes equations on surfaces.

In addition,
\begin{equation} \label{eq:matpartu}
\matpartu\, \zeta = \zeta_t + \vec u \cdot\nabla\,\zeta
\qquad \forall\ \zeta \in H^1(\GT)
\end{equation}
denotes the material time derivative of $\zeta$ on $\Gamma(t)$.
We compute $\matpartu\, \zeta$ with the help of an extension of $\zeta$ to a 
neighbourhood of $\GT$.
Here, we stress
that the derivative in (\ref{eq:matpartu}) is well-defined, and depends only on
the values of $\zeta$ on $\GT$, even though $\zeta_t$ and
$\nabla\,\zeta$ do not make sense separately for a function defined
on $\GT$; see e.g.\ \cite[p.\ 324]{DziukE13}.
The system \eqref{eq:1}, (\ref{eq:sigmaG}) is closed with the initial conditions
\begin{equation} \label{eq:1init}
\Gamma(0) = \Gamma_0 \,, \quad 
\rho\,\vec u(\cdot,0) = \rho\,\vec u_0 \quad \mbox{on } \Gamma_0
\,,
\end{equation}
where $\Gamma_0 \subset \bR^3$ and
$\vec u_0 : \Gamma_0 \to \bR^3$ are given initial data satisfying
$\rho\,\nabs\cdot\vec u_0 = 0$ on $\Gamma_0$.
Of course, in the case $\rho = 0$ the 
initial data $\vec u_0$ is not needed.

It is not difficult to show that (\ref{eq:1c}) leads to the conservation of
the total surface area $\mathcal{H}^{2}(\Gamma(t))$, see
Section~\ref{sec:weak} below for the relevant proofs. 
Furthermore, we note that
\begin{align}
\nabs\cdot\mat\sigma & =
2\,\mu\,\nabs\cdot\mat D_s (\vec u) -
\nabs\cdot[ p\,\mat{\mathcal{P}}] 
= 2\,\mu\,\nabs\cdot\mat D_s (\vec u) 
- \nabs\,p - \varkappa\,p\,\vec\nu\,.
\label{eq:nabssigma}
\end{align}
We recall that for
the mean curvature $\varkappa$ 
we use the sign
convention that $\varkappa$ is negative 
in case that $\Gamma(t)$ encloses a convex set.
In particular, it holds that
\begin{equation*} 
\Delta_s\, \vec\id = \varkappa\, \vec\nu 
\qquad \mbox{on $\Gamma(t)$}\,,
\end{equation*}
where $\Delta_s = \nabs\cdot\nabs$ is the Laplace--Beltrami operator on 
$\Gamma(t)$. 
We observe from \eqref{eq:nabssigma} that any term of the form
$\hat\alpha\varkappa\vec\nu$ on the right hand side of \eqref{eq:1b} would have no
effect on the solutions $\vec u$ and $\Gamma$, since if 
$(\vec u, p, \Gamma)$ is a solution to \eqref{eq:1} with the additional term
$\hat\alpha\varkappa\vec\nu$, then $(\vec u, p + \hat\alpha, \Gamma)$
solves \eqref{eq:1} as stated, i.e.\ without the extra term.

Assuming there are two solutions $\{(\Gamma(t), \vec u(\cdot,t),
p^{(i)}(\cdot,t))\}_{t \in [0,T]}$, $i=1,2$, 
to the problem
\eqref{eq:1}, (\ref{eq:sigmaG}) and
(\ref{eq:1init}), then it follows from (\ref{eq:nabssigma}) that
\begin{alignat}{2}
\nabs\,\bar p + \varkappa\,\bar p\,\vec\nu 
& = \vec 0
&&\qquad \text{on}\ \Gamma(t)\,,
\label{eq:barpb}
\end{alignat}
where
$\bar p =  p^{(1)} - p^{(2)}$.
Since $\nabs\,\bar p$ is tangential, we obtain that
$\nabs\,\bar p=\vec 0$, 
and hence 
$\bar p$ is a constant. Moreover, (\ref{eq:barpb}) implies that
$\varkappa\,\bar p = 0$. So unless
$\varkappa$ vanishes constantly, 
which cannot happen if $\Gamma(t)$ is compact,
then $\bar p = 0$. Hence $p$ is unique.
We remark that this is in contrast to the model
\eqref{eq:esns} with $\vec u \cdot \vec \nu = V$ a given 
velocity, where the pressure is only unique up to an additive constant.
See e.g.\ \cite{JankuhnORZ21,OlshanskiiRZ22} for more details.

Finally, we recall that the source term $\vec f_\Gamma = f_\Gamma\,\vec\nu$ 
in (\ref{eq:esns}), with $f_\Gamma$ defined in (\ref{eq:fGamma}), 
is the first variation of $E(\Gamma(t))$, i.e.\
\begin{equation*} 
\tfrac12\,\ddt \left\langle \varkappa, \varkappa \right\rangle_{\Gamma(t)}
= - \left\langle f_\Gamma , {\mathcal{V}} \right\rangle_{\Gamma(t)}
= - \left\langle \vec f_\Gamma , \vec{\mathcal{V}} \right\rangle_{\Gamma(t)}\,,
\end{equation*}
where $\langle \cdot, \cdot \rangle_{\Gamma(t)}$ denotes the 
$L^2$--inner product on $\Gamma(t)$.
Similarly to \cite{nsns}, in this paper we will make use of
the stable approximation of Willmore flow introduced in \cite{Dziuk08}, which
is based on a discretization of the curvature vector 
$\vec\varkappa = \varkappa\,\vec\nu$, and on the identity 
\begin{align}
\tfrac12\,\ddt \left\langle \vec \varkappa , \vec\varkappa 
 \right\rangle_{\Gamma(t)} & = 
-\left\langle \nabs\,\vec \varkappa , \nabs\,\vec{\mathcal{V}} 
 \right\rangle_{\Gamma(t)}
-\left\langle \nabs\cdot \vec \varkappa, \nabs\cdot \vec{\mathcal{V}} 
 \right\rangle_{\Gamma(t)}
- \tfrac{1}{2}\left\langle |\vec \varkappa|^2\,\nabs\, \vec\id, 
 \nabs\,\vec {\mathcal{V}} \right\rangle_{\Gamma(t)}
\nonumber \\ & \quad
+ 2 \left\langle (\nabs\,\vec \varkappa)^T, \mat D_s(\vec{{\mathcal{V}}})\,
 (\nabs\,\vec\id)^T \right\rangle_{\Gamma(t)} ,
\label{eq:dagger}
\end{align}
where we note that our notation is such that 
$\nabs\,\vec\chi = (\nabla_{\!\Gamma}\,\vec\chi)^T$, with 
$\nabla_{\!\Gamma}\,\vec\chi=\left( \partial_{s_i}\, \chi_j \right)_{i,j=1}^3$ 
defined as in \cite{Dziuk08}.

\begin{remark} \label{rem:general}
The following generalizations of \eqref{eq:1} are possible.
\begin{itemize}
\item
By allowing a term $\lambda \vec\nu$ on the right hand side of \eqref{eq:1b},
and choosing $\lambda(t) \in \bR$ such that 
$\int_{\Gamma(t)} \vec u\cdot\vec\nu \dH{2} = 0$, 
we can include conservation of the
volume enclosed by $\Gamma$ into the model. Then \eqref{eq:1} serves as a very 
simple model for fluidic biomembranes in the spirit of
\cite{ArroyoS09,nsns}, in particular in the limit of bulk densities and bulk
viscosities going to zero. 
We also refer to \cite{RahimiA12,RodriguesAMB15,KrauseV23},
where very similar models have been considered.
\item
Adding a term $\beta\vec u$ to the left hand side of \eqref{eq:1b} allows the
modelling of a Brinkman-type equation, see also \cite{OlshanskiiRZ21}. 
\item
Making $\rho$, $\mu$ and $\alpha$ phase-dependent
and introducing an interfacial energy on the surface, it is possible to consider
two-phase surface flows. See e.g.\ \cite{BachiniKNV23,ElliottS25} 
and also the related \cite{nsns2phase}.
\end{itemize}
\end{remark}

\begin{remark} \label{rem:killingfields}
So called Killing fields on an evolving surface $\Gamma(t)$ are nontrivial
solutions to the homogeneous tangential Stokes equations. That is, solutions of
the form $(\Gamma(t), \vec u(\cdot, t), p(t)) = (\Gamma_0, \vec u_0, 0)$, where
$\vec u_0 = \mat{\mathcal{P}}\,\vec u_0$ with $\mat D_s(\vec u_0) = \mat 0$.
For the mathematical analysis these pose intricate issues
for surface Stokes equations.
By way of example, on the unit sphere $\bS^2$ the space of Killing fields 
has dimension 3: one for each axis of rotation.
For a general smooth surface in $\bR^3$ it can be shown that the space of
Killing fields has dimension at most 3, see e.g.\ \cite{Sakai96}.
\end{remark}

\setcounter{equation}{0}
\section{Weak formulation} \label{sec:weak}

We define, similarly to (\ref{eq:matpartu}), the following time derivative that
follows the parameterization $\vec x(\cdot, t)$ of $\Gamma(t)$, rather than
$\vec u$. In particular, we let
\begin{equation*} 
\matpartx\, \zeta = \zeta_t + \vec{\mathcal{V}} \cdot\nabla\,\zeta
\qquad \forall\ \zeta \in H^1(\GT);
\end{equation*}
where we stress once again that this definition is well-defined, even though
$\zeta_t$ and $\nabla\,\zeta$ do not make sense separately for a
function $\zeta \in H^1(\GT)$.
On recalling (\ref{eq:matpartu}) we obtain that
$\matpartx = \matpartu$ if 
$\vec{\mathcal{V}} = \vec u$ 
on $\Gamma(t)$.
Moreover, for later use we note that
\begin{subequations}
\begin{equation} 
\langle \zeta, \nabs\cdot \vec\eta \rangle_{\Gamma(t)}
+ \langle \nabs\,\zeta , \vec\eta \rangle_{\Gamma(t)}
= - \langle \zeta\,\vec\eta , \vec\varkappa \rangle_{\Gamma(t)} 
\qquad \forall\ \zeta \in H^1(\Gamma(t)),\, \vec \eta \in [H^1(\Gamma(t))]^3
\end{equation}
and
\begin{equation} \label{eq:DElem5.2}
\ddt \left\langle \chi, \zeta \right\rangle_{\Gamma(t)}
 = \left\langle \matpartx\,\chi, \zeta \right\rangle_{\Gamma(t)}
 + \left\langle \chi, \matpartx\,\zeta \right\rangle_{\Gamma(t)}
+ \left\langle \chi\,\zeta, \nabs\cdot\vec{\mathcal{V}} 
 \right\rangle_{\Gamma(t)}
\qquad \forall\ \chi,\zeta \in H^1(\GT)\,,
\end{equation}
\end{subequations}
see Definition~2.11 and Lemma~5.2 in \cite{DziukE13}, respectively.

The most natural weak formulation of the system \eqref{eq:esns} 
uses the fluidic tangential velocity for the evolution of $\Gamma(t)$, and so
the third equation in \eqref{eq:esns} is replaced by $\vec {\mathcal{V}} = \vec u$ on $\Gamma(t)$.
It then follows from \eqref{eq:esns} that
\begin{align*}
& \rho
\left\langle \matpartu\, \vec u, \vec\xi 
\right\rangle_{\Gamma(t)}
+ 2\,\mu \left\langle \mat D_s (\vec u) , \mat D_s (\vec \xi)
\right\rangle_{\Gamma(t)}
- \left\langle p, \nabs\cdot\vec\xi \right\rangle_{\Gamma(t)}
= \left\langle \vec g + \alpha \vec f_\Gamma, \vec \xi \right\rangle_{\Gamma(t)}
\nonumber \\ & \hspace{11cm}
\quad \forall\ \vec \xi \in [H^1(\GT)]^3\,,
\end{align*}
where we have noted for symmetric matrices $\mat A \in \bR^{3 \times 3}$
that $\mat{\mathcal{P}}\,\mat A\,\mat{\mathcal{P}} : \mat B = 
\mat{\mathcal{P}}\,\mat A \, \mat{\mathcal{P}} : \tfrac12\,\mat{\mathcal{P}}\,(\mat B + \mat B^T)\,
\mat{\mathcal{P}}$ for all $\mat B \in \bR^{3 \times 3}$.
The resulting weak formulation of \eqref{eq:esns} 
is then given as follows.
Find $\Gamma(t) = \vec x(\Upsilon, t)$ for $t\in[0,T]$ 
with velocity $\vec{\mathcal{V}} \in [L^2(\GT)]^3$
and $\vec{\mathcal{V}}(\cdot, t) \in [H^1(\Gamma(t)]^3$ for almost all 
$t \in (0,T)$, and functions 
$\vec u \in [H^1(\GT)]^3$, 
$p \in L^2(\GT)$,
$\vec\varkappa \in [H^1(\GT)]^3$ and
$\vec f_\Gamma \in [L^2(\GT)]^3$ such that 
for almost all $t \in (0,T)$ it holds that
\begin{subequations} \label{eq:weakGD}
\begin{align}
& \rho \left\langle \matpartx\,\vec u , \vec \xi 
 \right\rangle_{\Gamma(t)}
+ 2\,\mu \left\langle \mat D_s (\vec u) , \mat D_s (\vec \xi)
\right\rangle_{\Gamma(t)}
- \left\langle p, \nabs\cdot\vec \xi \right\rangle_{\Gamma(t)}
 = \left\langle \vec g + \alpha \vec f_\Gamma, \vec \xi \right\rangle_{\Gamma(t)} \nonumber \\ & \hspace{9cm}
\quad \forall\ \vec \xi \in [H^1(\GT)]^3\,,
\label{eq:weakGDa} \\
& \left\langle \nabs\cdot\vec u, \eta \right\rangle_{\Gamma(t)}
 = 0  \quad \forall\ \eta \in L^2(\Gamma(t))\,, \label{eq:weakGDc} \\
&  \left\langle \vec{\mathcal{V}} - \vec u, \vec\chi \right\rangle_{\Gamma(t)} = 0
 \quad\forall\ \vec\chi \in [L^2(\Gamma(t))]^3\,,
\label{eq:weakGDd} \\
& \left\langle \vec\varkappa, \vec\eta \right\rangle_{\Gamma(t)}
+ \left\langle \nabs\,\vec\id, \nabs\,\vec \eta \right\rangle_{\Gamma(t)}
 = 0  \quad\forall\ \vec\eta \in [H^1(\Gamma(t))]^3\,, \label{eq:weakGDe} \\
& \left\langle \vec f_\Gamma, \vec \chi \right\rangle_{\Gamma(t)} =
\left\langle \nabs\,\vec \varkappa , \nabs\,\vec \chi \right\rangle_{\Gamma(t)}
+\left\langle \nabs\cdot \vec \varkappa, \nabs\cdot \vec \chi 
\right\rangle_{\Gamma(t)}
+ \tfrac{1}{2}\left\langle |\vec \varkappa|^2\,\nabs\, \vec\id, 
\nabs\,\vec \chi \right\rangle_{\Gamma(t)}
\nonumber \\ & \hspace{4cm}
- 2 \left\langle (\nabs\,\vec \varkappa)^T, \mat D_s(\vec{\chi})\,
 (\nabs\,\vec\id)^T \right\rangle_{\Gamma(t)}
\quad\forall\ \vec \chi \in [H^1(\Gamma(t))]^3 \,, \label{eq:weakGDf}
\end{align}
\end{subequations}
as well as initial conditions for $\Gamma$ and $\vec u$.

In what follows we would like to derive an energy bound for a solution of
\eqref{eq:weakGD}. 
All of the following considerations are formal, in the
sense that we make the appropriate assumptions about the existence, 
boundedness and regularity of a solution to \eqref{eq:weakGD}. 
Firstly, it follows from (\ref{eq:DElem5.2}), (\ref{eq:weakGDd}) and 
(\ref{eq:weakGDc}) with $\eta = |\vec u|^2$ that
\begin{align} \label{eq:dtuu}
\tfrac12\,\rho \,\ddt
\left\langle \vec u, \vec u \right\rangle_{\Gamma(t)}
& = \tfrac12\,\rho \left\langle \matpartx\,|\vec u|^2,1
\right\rangle_{\Gamma(t)}
+ \tfrac12\,\rho \left\langle \nabs\cdot\vec{\mathcal{V}}, |\vec u|^2
\right\rangle_{\Gamma(t)}
\nonumber \\ &
 = \rho \left\langle \matpartx\,\vec u, \vec u
\right\rangle_{\Gamma(t)}
+ \tfrac12\,\rho \left\langle \nabs\cdot \vec u, 
|\vec u|^2 \right\rangle_{\Gamma(t)}
= \rho \left\langle \matpartx\,\vec u, \vec u
\right\rangle_{\Gamma(t)} .
\end{align}
Now choosing $\vec\xi = \vec u$ in (\ref{eq:weakGDa})
and $\eta = p(\cdot, t)$ in (\ref{eq:weakGDc}) yields,
on combining with (\ref{eq:dtuu}), that
\begin{align} \label{eq:d1}
& \tfrac12\,\ddt \left( 
\rho \left\langle \vec u, \vec u \right\rangle_{\Gamma(t)} \right)
+ 2\,\mu\left\langle \mat D_s(\vec u), \mat D_s(\vec u) 
\right\rangle_{\Gamma(t)} 
= \left\langle\vec g, \vec u\right\rangle_{\Gamma(t)} 
+ \alpha \left\langle \vec f_\Gamma , \vec u \right\rangle_{\Gamma(t)}.
\end{align}
Combining (\ref{eq:d1}) with (\ref{eq:dagger}), on choosing
$\vec\chi = \vec f$ in (\ref{eq:weakGDd}) and
$\vec\chi = \vec{\mathcal{V}}$ in (\ref{eq:weakGDf}), yields that
\begin{align*} 
& \tfrac12\,\ddt \left( 
\rho \left\langle \vec u, \vec u \right\rangle_{\Gamma(t)} 
+ \alpha \left\langle \vec\varkappa, \vec\varkappa \right\rangle_{\Gamma(t)} 
 \right)
+ 2\,\mu\left\langle \mat D_s(\vec u), \mat D_s(\vec u) 
\right\rangle_{\Gamma(t)} 
= \left\langle\vec g, \vec u\right\rangle_{\Gamma(t)} \,.
\end{align*}

Moreover, we note that it
immediately follows from (\ref{eq:DElem5.2}) and (\ref{eq:weakGDc},d) that
\begin{equation} \label{eq:areacons}
\ddt\, \mathcal{H}^{2} (\Gamma(t)) = \ddt\,\langle 1, 1 \rangle_{\Gamma(t)}
= \left\langle 1, \nabs\cdot\vec{\mathcal{V}} \right\rangle_{\Gamma(t)}
= \left\langle 1, \nabs\cdot\vec u \right\rangle_{\Gamma(t)}
= 0 \,.
\end{equation}

We note the following LBB condition.
\begin{lemma} \label{lem:LBB}
Let $(\Gamma(t))_{t\in [0,T]}$ be sufficiently smooth.
Let $\mathbb V(\Gamma(t)) = \{ \vec\xi \in [L^2(\Gamma(t)]^3 :
\mat{\mathcal{P}}\,\vec\xi \in [H^1(\Gamma(t)]^3\}$,  
and let $\|\cdot\|_{0,\Gamma(t)}$
and $\|\cdot\|_{1,\Gamma(t)}$ denote the usual Sobolev norms on $\Gamma(t)$. 
Then there exists a constant $C>0$ such that
\begin{equation} \label{eq:LBB}
\inf_{\eta \in L^2(\Gamma(t))} 
\sup_{\vec \xi \in \mathbb V(\Gamma(t))}
\frac{ \left\langle \eta, \nabs \cdot\vec \xi \right\rangle_{\Gamma(t)}}
{\| \eta \|_{0,\Gamma(t)}\, (\| \vec\xi \|_{0,\Gamma(t)} +
\|\mat{\mathcal{P}}\,\vec \xi\|_{1,\Gamma(t)})} 
\geq C.
\end{equation}
\end{lemma}
\begin{proof}
Let an arbitrary $\eta \in L^2(\Gamma(t))$ be given.
We first choose $\xi_n \in L^2(\Gamma(t))$ such that
\begin{equation} \label{eq:xin}
\eta + \xi_n\,\varkappa \in L^2_0(\Gamma(t)) := \{
\zeta \in L^2(\Gamma(t)) : \left\langle \zeta,1 \right\rangle_{\Gamma(t)} 
= 0 \}.
\end{equation}
This is always possible, 
due to $\Gamma(t)$ being a closed surface without boundary. In fact,
due to our sign convention, $\varkappa$ needs to be negative on a set of
positive measure. We define $E = \{ \vec z \in \Gamma(t) : \varkappa(\vec z) <
0\}$, and from now on set 
$\xi_n = - \frac{\left\langle \eta, 1 \right\rangle_{\Gamma(t)}}{
\left\langle \varkappa, \charfcn{E} \right\rangle_{\Gamma(t)}} \charfcn{E}$,
which clearly satisfies \eqref{eq:xin}. For later use we note that with this
choice
\begin{equation} \label{eq:xinbound}
\|\xi_n\|_{0,\Gamma(t)} \leq C_1 \|\eta\|_{0,\Gamma(t)}
\end{equation}
for a suitable constant $C_1 >0$ independent of $\eta$.
Next we let $w \in H^1(\Gamma(t)) \cap L^2_0(\Gamma(t))$ 
be the solution to $-\Delta_s\, w = \eta +
\xi_n\,\varkappa$. The existence and uniqueness of $w$ immediately follows from 
\eqref{eq:xin} and the Lax--Milgram theorem, together with the Poincar\'e
inequality. Elliptic regularity theory implies that
$w \in H^2(\Gamma(t))$ with
\begin{equation} \label{eq:ellreg}
\|w\|_{2,\Gamma(t)} \leq C_2 \|\Delta_s\,w\|_{0,\Gamma(t)}
= C_2 \|\eta + \xi_n\,\varkappa\|_{0,\Gamma(t)}.
\end{equation}
Then we set
\begin{equation*} 
\vec\xi = \xi_n\,\vec\nu + \vec\xi_\subT\,,\quad
\text{with } \vec\xi_\subT = - \nabs\,w
\end{equation*}
and 
compute from \eqref{eq:nabs2} that 
$\nabs\cdot\vec\xi = -\xi_n\,\varkappa + \nabs\cdot\vec\xi_\subT
= -\xi_n\,\varkappa - \Delta_s\, w = \eta$. Moreover, we have from
\eqref{eq:ellreg} and \eqref{eq:xinbound} that 
\begin{equation*} 
\|\vec\xi_\subT\|_{1,\Gamma(t)} 
\leq C_2 \|\eta + \xi_n\,\varkappa\|_{0,\Gamma(t)}
\leq (C_2 + C_1 \sup_{\Gamma(t)} |\varkappa|) \|\eta \|_{0,\Gamma(t)}
\leq C_3 \|\eta \|_{0,\Gamma(t)},
\end{equation*}
and note furthermore that
$\| \vec\xi \|_{0,\Gamma(t)} \leq \| \xi_n \|_{0,\Gamma(t)} + 
\| \vec\xi_\subT \|_{0,\Gamma(t)} \leq (C_1+C_3) \| \eta \|_{0,\Gamma(t)}$.
Consequently, it holds that
\begin{align*}
\frac{ \left\langle \eta, \nabs \cdot\vec \xi \right\rangle_{\Gamma(t)}}
{\| \eta \|_{0,\Gamma(t)}\, (\| \vec\xi \|_{0,\Gamma(t)} +
\|\mat{\mathcal{P}}\,\vec \xi\|_{1,\Gamma(t)})} 
& = \frac{ \| \eta \|_{0,\Gamma(t)}}
{\| \vec\xi \|_{0,\Gamma(t)} +
\|\vec \xi_\subT\|_{1,\Gamma(t)}} \nonumber \\ &
\geq
\frac{ \| \eta \|_{0,\Gamma(t)}}
{(C_1 + C_3) \| \eta \|_{0,\Gamma(t)}
+ C_3 \|\eta\|_{0,\Gamma(t)}}
\geq (C_1 + 2 C_3)^{-1}\,.
\end{align*}
\end{proof}

\setcounter{equation}{0}
\section{Semidiscrete finite element approximation} \label{sec:sd}

We begin by extending the notion of an {\em evolving polyhedral
hypersurface} from \cite[Definition~68]{bgnreview} to that of an
{\em evolving discrete $P\ell$ hypersurface}, for $\ell \geq 1$.
To this end, let
$(\widehat\Gamma^h(t))_{[0,T]}$ be an evolving polyhedral
hypersurface in $\bR^3$, such that 
$\widehat\Gamma^h(t)=\bigcup_{j=1}^{J} \overline{\widehat\sigma^h_j(t)}$, where 
$\{\widehat\sigma^h_j(t)\}_{j=1}^{J}$ is a
family of mutually disjoint open triangles with vertices
$\{\vec{\mathfrak q}^h_k(t)\}_{k=1}^{\widehat K}$.
Then for each fixed $t\in [0,T]$, we let
$\Gamma^h(t)=\bigcup_{j=1}^{J} \overline{\sigma^h_j(t)}$, where there
exists a homeomorphism $\mathcal F: \widehat\Gamma^h(t) \to \Gamma^h(t)$
such that $\mathcal F_{\mid_{\widehat\sigma^h_j(t)}}$ is a polynomial of degree
$\ell$,
$\mathcal F(\widehat\sigma^h_j(t)) 
= \sigma^h_j(t)$, $j=1,\ldots, J$, and
$\mathcal F(\vec{\mathfrak q}^h_k(t)) = \vec{\mathfrak q}^h_k(t)$,
$k=1,\ldots,\widehat K$.
We denote by $\{\vec p^h_k(t)\}_{k=1}^{K}$ the set of
Lagrange nodes on $\widehat\Gamma^h(t)$ for the space of continuous finite
elements of degree $\ell$. That is, it contains for each triangle the vertices,
on each edge $\ell-1$ equally spaced points, and 
$\binom{\ell+2}\ell - 3\ell = \frac12(\ell-1)(\ell-2)$
interior points. We then define the Lagrange nodes of the curved
surface $\Gamma^h(t)$ by $\{\vec q^h_k(t)\}_{k=1}^{K}$, where
$\vec{q}^h_k(t) = \mathcal F(\vec p^h_k(t))$,
$k=1,\ldots, K$.
The final assumption for $(\Gamma^h(t))_{[0,T]}$ to be an evolving 
discrete $P\ell$ hypersurface is then that 
the functions $t \mapsto \vec q_k(t)$, are elements of $C^1([0,T])$ for
$k=1,\ldots,K$.

We then define the following finite element spaces. Let
\begin{equation} \label{eq:Shr}
 S^h_r(\Gamma^h(t)) = \{\chi \in C^0(\Gamma^h(t)) : 
\chi_{\mid_{\sigma^h_j(t)}} \circ \mathcal F_{\mid_{\widehat\sigma^h_j(t)}}
\in \mathbb{P}^r(\widehat\sigma^h_j(t)) \quad \forall\ j=1,\ldots, J\} 
\,, 
\end{equation}
for $r=1,\ldots,\ell$.
Let $\{\chi^h_k(\cdot,t)\}_{k=1}^{K}$ 
denote the standard basis of $S^h_\ell(\Gamma^h(t))$, i.e.\
\begin{equation} \label{eq:bf}
\chi^h_k(\vec q^h_l(t),t) = \delta_{kl}\qquad
\forall\ k,l \in \{1,\ldots,K\}\,,\ t \in [0,T]\,.
\end{equation}
We define the discrete material velocity for $\vec z \in \Gamma^h(t)$ by
\begin{equation} \label{eq:Xht}
\vec{\mathcal{V}}^h(\vec z, t) = \sum_{k=1}^{K}
\left[\ddt\,\vec q^h_k(t)\right] \chi^h_k(\vec z, t) \,,
\end{equation}
and the introduce the material time derivative
\begin{equation*} 
\matpartxh\, \zeta = \zeta_t + \vec{\mathcal{V}}^h\cdot\nabla\,\zeta
\qquad\forall\ \zeta \in H^1(\GhT)\,,
\quad\text{where}\quad
\GhT = \bigcup_{t \in [0,T]} \Gamma^h(t) \times \{t\}\,.
\end{equation*}

For later use, for $r=1,\ldots,\ell$, we also introduce the finite element 
spaces
\begin{align*}
S^h_r(\GhT) & = \{ \chi \in C^0(\GhT) : 
\chi(\cdot, t) \in S^h_r(\Gamma^h(t)) \quad \forall\ t \in [0,T] \}\,, \\
S^h_{r,T}(\GhT) & = \{ \chi \in S^h_r(\GhT) : \matpartxh\,\chi \in C^0(\GhT) \}\,.
\end{align*}

On differentiating (\ref{eq:bf}) with respect to $t$, it immediately follows
that
\begin{equation} \label{eq:mpbf}
\matpartxh\, \chi^h_k = 0
\quad\forall\ k \in \{1,\ldots,K\}\,.
\end{equation}
It follows directly from (\ref{eq:mpbf}) that
\begin{equation*} 
\matpartxh\,\zeta(\cdot,t) = \sum_{k=1}^{K} \chi^h_k(\cdot,t)\,
\ddt\,\zeta_k(t) \quad \text{on}\ \Gamma^h(t)
\end{equation*}
for $\zeta(\cdot,t) = \sum_{k=1}^{K} \zeta_k(t)\,\chi^h_k(\cdot,t)
\in S^h_\ell(\Gamma^h(t))$, and hence
$\matpartxh\,\vec\id = \vec{\mathcal{V}}^h$ on $\Gamma^h(t)$.
We refer to \cite{Heine03PhD,BrennerS02,Demlow09,ElliottR21} on more details
on higher order (evolving) parametric finite element spaces.
In particular, we recall from \cite[Lem.~9.9]{ElliottR21} that
\begin{equation} \label{eq:ERlem99}
\ddt \langle \eta, 1 \rangle_{\Gamma^h(t)}
 = \langle \matpartxh\,\eta, 1 \rangle_{\Gamma^h(t)}
+ \langle \eta, \nabs\cdot\vec{\mathcal{V}}^h \rangle_{\Gamma^h(t)}
\qquad \forall\ \eta \in C^1(\GhT)\,.
\end{equation}

We introduce, similarly to \eqref{eq:diff2}, 
\begin{subequations} \label{eq:diffh}
\begin{equation} \label{eq:Psh}
\mat{\mathcal{P}}^h = \mat\Id - \vec\nu^h \otimes \vec\nu^h
\quad\text{on}\ \Gamma^h(t)\,,
\end{equation}
where $\vec{\nu}^h(t) \in [L^\infty(\Gamma^h(t))]^3$ denotes the unit normal
to $\Gamma^h(t)$, and
\begin{equation} \label{eq:Dsh} 
\mat D_s^h(\vec \eta) = \tfrac12\,\mat{\mathcal{P}}^h\,
(\nabs\,\vec\eta + (\nabs\,\vec \eta)^T)\,\mat{\mathcal{P}}^h
\quad\text{on}\ \Gamma^h(t)\,,
\end{equation}
\end{subequations}
where here $\nabs = \mat{\mathcal{P}}^h \,\nabla$ 
denotes the surface gradient on $\Gamma^h(t)$.

In what follows we will introduce a finite element approximation
for (\ref{eq:esns}).
When designing such a finite element approximation, a
careful decision has to be made about the {\em discrete tangential velocity} of
$\Gamma^h(t)$. The most natural choice is to select the velocity of the fluid,
i.e.\ (\ref{eq:weakGDd}) is appropriately discretized, and that is the approach
we adopt in this paper.
In particular, for the discrete approximation of the velocity, $\vec U^h$,
we will employ the $P\ell$-isoparametric elements $[S^h_{\ell}(\Gamma^h(t))]^3$
defined in \eqref{eq:Shr}.
As a consequence, and with a view towards well-posedness considerations on the
fully discrete level, we will seek for the surface pressure in the space
$S^h_{\ell-1}(\Gamma^h(t))$, so that the pressure-velocity pair corresponds
to the classical $P\ell$--$P(\ell-1)$ Taylor--Hood element in the flat case.

Overall, we then obtain the following semidiscrete
continuous-in-time finite element approximation, which is the semidiscrete
analogue of the weak formulation \eqref{eq:weakGD}.

Let $\ell \geq 2$ and choose a fixed parameter $\theta\in\{0,1\}$.
Given $\Gamma^h(0)$ and $\vec U^h(\cdot,0) \in [S^h_{\ell}(\Gamma^h(0))]^3$,
find an evolving discrete $P\ell$ hypersurface
$(\Gamma^h(t))_{t\in[0,T]}$ such that $\vec \id_{\mid_{\Gamma^h(t)}} \in \qVht$
for $t \in [0,T]$, and functions $\vec U^h \in [S^h_{\ell,T}(\GhT)]^3$, 
$P^h \in S^h_{\ell-1}(\GhT)$, $\vec\kappa^h \in [S^h_\ell(\GhT)]^3$ and
$\vec F^h \in [S^h_\ell(\GhT)]^3$
such that for almost all $t \in (0,T)$ it holds that
\begin{subequations} \label{eq:qsdGD}
\begin{align}
& \rho \left\langle \matpartxh\,\vec U^h, \vec \xi \right\rangle_{\Gamma^h(t)}
+ \tfrac\theta2 \rho \left\langle \nabs\cdot\vec U^h, \vec U^h \cdot \vec\xi 
\right\rangle_{\Gamma^h(t)} 
+ 2\,\mu \left\langle \mat D_s^h \vec U^h , 
\mat D_s^h \vec \xi \right\rangle_{\Gamma^h(t)}
- \left\langle P^h , 
\nabs\cdot\vec \xi \right\rangle_{\Gamma^h(t)}
\nonumber \\ & \hspace{4cm}
= \left\langle\vec g^h,\vec\xi \right\rangle_{\Gamma^h(t)} 
 + \alpha \left\langle \vec F^h,\vec\xi \right\rangle_{\Gamma^h(t)}
\qquad \forall\ \vec\xi \in [S^h_{\ell,T}(\GhT)]^3\,, \label{eq:qsdGDa}\\
& \left\langle \nabs\cdot\vec U^h, \eta 
\right\rangle_{\Gamma^h(t)}  = 0 
\qquad \forall\ \eta \in \Wht\,,
\label{eq:qsdGDc} \\
& \left\langle \vec{\mathcal{V}}^h ,
\vec\chi \right\rangle_{\Gamma^h(t)}
= \left\langle \vec U^h, \vec\chi \right\rangle_{\Gamma^h(t)}
 \qquad\forall\ \vec\chi \in \qVht\,,
\label{eq:qsdGDd} \\
& \left\langle \vec\kappa^h , \vec\eta \right\rangle_{\Gamma^h(t)}
+ \left\langle \nabs\,\vec\id, \nabs\,\vec \eta \right\rangle_{\Gamma^h(t)}
 = 0  \qquad\forall\ \vec\eta \in \qVht\,,\label{eq:qsdGDe} \\
& \left\langle \vec F^h, \vec\chi \right\rangle_{\Gamma^h(t)} =
\left\langle \nabs\,\vec{\kappa}^h , \nabs\,\vec\chi \right\rangle_{\Gamma^h(t)}
+ \left\langle \nabs\cdot\vec{\kappa}^h , \nabs\cdot\vec\chi \right\rangle_{\Gamma^h(t)}
+\tfrac12 \left\langle |\vec{\kappa}^h|^2\,\nabs\,\vec\id, 
 \nabs\,\vec\chi \right\rangle_{\Gamma^h(t)}
 \nonumber \\ & \hspace{3cm}
-2 \left\langle (\nabs\,\vec\kappa^h)^T , \mat D_s^h(\vec\chi)\,
 (\nabs\,\vec\id)^T \right\rangle_{\Gamma^h(t)}
\qquad \forall\ \vec\chi \in \qVht\,, \label{eq:qsdGDf} 
\end{align}
\end{subequations}
where we recall \eqref{eq:Xht}.
We observe that while \eqref{eq:qsdGD} with $\theta=0$ is the most natural
semidiscrete approximation of \eqref{eq:weakGD}, also \eqref{eq:qsdGD} with
$\theta=1$ is a consistent approximation of
\eqref{eq:weakGD} due to the corresponding term vanishing in the continuous
setting. And it will turn out that only with the choice $\theta=1$ we are able 
to prove a general energy stability result for \eqref{eq:qsdGD}. 

\begin{theorem} \label{thm:stabqGD}
Let $\{(\Gamma^h, \vec U^h, P^h, \vec\kappa^h, \vec F^h)(t)
\}_{t\in[0,T]}$ 
be a solution to \eqref{eq:qsdGD}. Then, 
if $\rho\,(\theta-1) = 0$, it holds that
\begin{align} \label{eq:qthmGD}
& \tfrac12\,\ddt \left( 
\rho \left\langle \vec U^h, \vec U^h \right\rangle_{\Gamma^h(t)}
+ \alpha \left\langle \vec\kappa^h,
\vec\kappa^h\right\rangle_{\Gamma^h(t)}
 \right)
+ 2\,\mu \left\langle \mat D_s^h \vec U^h , 
\mat D_s^h \vec U^h \right\rangle_{\Gamma^h(t)}
= \left\langle\vec g, \vec U^h\right\rangle_{\Gamma^h(t)}\,.
\end{align}
It also holds that
\begin{equation} \label{eq:qareaconsh}
\ddt\, \mathcal{H}^{2} (\Gamma^h(t)) = 0 \,.
\end{equation}
\end{theorem} 
\begin{proof}
Choosing
$\vec \xi = \vec U^h$ in (\ref{eq:qsdGDa})
and $\eta = P^h$ in (\ref{eq:qsdGDc}) yields that
\begin{align} \label{eq:qlemGD}
& \rho \left\langle \matpartxh\,\vec U^h, \vec U^h
\right\rangle_{\Gamma^h(t)}
+ \tfrac\theta2 \rho \left\langle \nabs\cdot\vec U^h, |\vec U^h|^2 
\right\rangle_{\Gamma^h(t)} 
+ 2\,\mu \left\langle \mat D_s^h \vec U^h , 
\mat D_s^h \vec U^h \right\rangle_{\Gamma^h(t)} \nonumber \\ & \qquad
= \left\langle \vec g^h,\vec U^h \right\rangle_{\Gamma^h(t)} + 
\alpha \left\langle \vec F^h,\vec U^h \right\rangle_{\Gamma^h(t)}.
\end{align}
Moreover, it is possible to show, similarly to (\ref{eq:dagger}), that
\begin{align*}
\tfrac12\,\ddt 
 \left\langle \vec\kappa^h, \vec\kappa^h\right\rangle_{\Gamma^h(t)} & = 
- \left\langle \nabs\,\vec{\kappa}^h , \nabs\,\vec{\mathcal{V}}^h 
 \right\rangle_{\Gamma^h(t)}
- \left\langle \nabs\cdot\vec{\kappa}^h , \nabs\cdot\vec{\mathcal{V}}^h 
 \right\rangle_{\Gamma^h(t)}
 \nonumber \\ & \quad
-\tfrac12 \left\langle |\vec{\kappa}^h|^2\,\nabs\,\vec\id, 
 \nabs\,\vec{\mathcal{V}}^h \right\rangle_{\Gamma^h(t)}
+ 2 \left\langle (\nabs\,\vec\kappa^h)^T, \mat D_s^h(\vec{\mathcal{V}}^h)
 \,(\nabs\,\vec\id)^T \right\rangle_{\Gamma^h(t)} , 
\end{align*}
see \cite{Dziuk08,ElliottS10}. Hence choosing $\vec\chi = \vec F^h$ in
(\ref{eq:qsdGDd}) and $\vec\chi = \vec{\mathcal{V}}^h$ in (\ref{eq:qsdGDf})
yields that
\begin{equation} \label{eq:qDziuk08f}
 \left\langle\vec F^h,\vec U^h \right\rangle_{\Gamma^h(t)} =
 \left\langle\vec F^h,\vec{\mathcal{V}}^h 
 \right\rangle_{\Gamma^h(t)}
= -\tfrac12\,\ddt 
 \left\langle \vec\kappa^h, \vec\kappa^h\right\rangle_{\Gamma^h(t)}.
\end{equation}
If $\rho = 0$, then the desired result (\ref{eq:qthmGD}) directly 
follows from combining (\ref{eq:qlemGD}) and (\ref{eq:qDziuk08f}). 
If $\rho > 0$, on the other hand, we note, similarly to (\ref{eq:dtuu}),
that (\ref{eq:ERlem99}) and \eqref{eq:qsdGDd} imply that
\begin{align} \label{eq:qdtUU}
\tfrac12\,\rho \,\ddt
\left\langle \vec U^h, \vec U^h \right\rangle_{\Gamma^h(t)}
& = \tfrac12\,\rho \left\langle \matpartxh\,[|\vec U^h|^2],1
\right\rangle_{\Gamma^h(t)}
+ \tfrac12\,\rho \left\langle \nabs\cdot\vec{\mathcal{V}}^h, |\vec U^h|^2
\right\rangle_{\Gamma^h(t)}
\nonumber \\ &
 = \rho \left\langle \matpartxh\,\vec U^h, \vec U^h
\right\rangle_{\Gamma^h(t)}
+ \tfrac12\,\rho \left\langle \nabs\cdot \vec U^h, 
|\vec U^h|^2 \right\rangle_{\Gamma^h(t)} .
\end{align}
Combining (\ref{eq:qlemGD}), (\ref{eq:qdtUU}) and (\ref{eq:qDziuk08f}) 
yields the desired result (\ref{eq:qthmGD}) for $\rho > 0$,
on recalling that the assumptions then mean that $\theta=1$.

It immediately follows from (\ref{eq:ERlem99}) that
\begin{equation*} 
\ddt\left\langle 1, 1 \right\rangle_{\Gamma^h(t)} = 
\left\langle 1, \nabs\cdot\vec{\mathcal{V}}^h \right\rangle_{\Gamma^h(t)}
=
\left\langle 1, \nabs\cdot \vec U^h \right\rangle_{\Gamma^h(t)} = 0 \,,
\end{equation*}
which proves the desired result (\ref{eq:qareaconsh}). 
\end{proof}

\begin{remark} \label{rem:theta0}
We observe that in the case $\theta=0$
it does not seem possible to mimic the stability proof from the continuous
level if $\rho>0$, recall \eqref{eq:dtuu}. The reason is that
$\eta = |\vec U^h|^2$ is usually not an admissible discrete 
pressure test function.
In practice, for the fully discrete variant introduced in Section~\ref{sec:fd},
the differences between simulations for $\theta=0$ and $\theta=1$ are 
negligible. We therefore believe that the restriction to $\theta=1$ for the
stability proof is a theoretical one, rather than an inherent instability of
the scheme with $\theta=0$.
\end{remark}

\setcounter{equation}{0}
\section{Fully discrete finite element approximation} \label{sec:fd}

In this section we consider a fully discrete variant of the scheme
\eqref{eq:qsdGD} from \S\ref{sec:sd}. Here we will
choose the time discretization such that existence and uniqueness of the
discrete solutions can be guaranteed, and such that we inherit as much of the
structure of the stable schemes in \cite{spurious,fluidfbp} as possible, see
also \cite{nsns} for more details.

We consider the partitioning $t_m = m\,\tau$, $m = 0,\ldots, M$, 
of $[0,T]$ into uniform time steps $\tau = T / M$.
Then let $\Gamma^m$ be a fully discrete approximation of $\Gamma^h(t_m)$,
with all the previously introduced notations for $\Gamma^h(t_m)$ naturally
carrying across to $\Gamma^m$. For example,
we have that $\Gamma^m=\bigcup_{j=1}^{J} \overline{\sigma^m_j}$, and
there exists a homeomorphism $\mathcal F: \widehat\Gamma^m \to \Gamma^m$
that is a piecewise polynomial of degree at most $\ell$,
with $\widehat\Gamma^m=\bigcup_{j=1}^{J} \overline{\widehat\sigma^m_j}$. 
The set of Lagrange nodes of $\Gamma^m$ is denoted by 
$\{\vec q^m_k\}_{k=1}^{K}$. 
Let
\begin{equation*} 
 S^h_r(\Gamma^m) = \{\chi \in C^0(\Gamma^m) : 
\chi_{\mid_{\sigma^m_j}} \circ \mathcal F_{\mid_{\widehat\sigma^m_j}}
\in \mathbb{P}^r(\widehat\sigma^m_j) \quad \forall\ j=1,\ldots, J\} 
\,, \quad r=1,\ldots,\ell\,.
\end{equation*}
Let $\{\chi^m_k\}_{k=1}^{K}$ 
denote the standard basis of $S^h_\ell(\Gamma^m)$, i.e.\
\begin{equation*} 
\chi^m_k(\vec q^m_l) = \delta_{kl}\qquad
\forall\ k,l \in \{1,\ldots,K\}\,.
\end{equation*}
We also introduce the standard interpolation operator
$\vec\pi^m: [C^0(\Gamma^m)]^3 \to [S^h_\ell(\Gamma^m)]^3$, as well as the
following pushforward operator for the discrete interfaces
$\Gamma^m$ and $\Gamma^{m-1}$, for $m=0,\ldots,M$. Here we set
$\Gamma^{-1}=\Gamma^0$. Let $\vec\Pi_{m-1}^{m} : [C^0(\Gamma^{m-1})]^3 \to
\qVh$ such that
\begin{equation*} 
\vec\Pi_{m-1}^{m}\,\vec z = \sum_{k=1}^K \chi^m_k\,\vec z(\vec q^{m-1}_k)\,,
\qquad \forall\ \vec z \in [C^0(\Gamma^{m-1})]^3\,,
\end{equation*}
for $m=1,\ldots,M$, and set $\vec\Pi_{-1}^{0} = \vec\pi^{0}$.

We also introduce the $L^2$--inner 
product $\langle\cdot,\cdot\rangle_{\Gamma^m}$ over
the current polyhedral surface $\Gamma^m$.
We introduce, similarly to \eqref{eq:diffh},
\begin{equation*} 
\mat{\mathcal{P}}^m = \mat\Id - \vec \nu^m \otimes \vec \nu^m
\quad\text{on}\ \Gamma^m\,,
\end{equation*}
and
\begin{equation*} 
\mat D_s^m(\vec \eta) = \tfrac12\,\mat{\mathcal{P}}^m\,
(\nabs\,\vec\eta + (\nabs\,\vec \eta)^T)\,\mat{\mathcal{P}}^m
\quad\text{on}\ \Gamma^m\,,
\end{equation*}
where here $\nabs = \mat{\mathcal{P}}^m \,\nabla$ 
denotes the surface gradient on $\Gamma^m$.

Throughout this paper, we will parameterize the new closed surface 
$\Gamma^{m+1}$ over $\Gamma^m$, with the help of a parameterization
$\vec X^{m+1} \in \qVh$, i.e.\ $\Gamma^{m+1} = \vec X^{m+1}(\Gamma^m)$.

Let $\ell \geq 2$ and choose a fixed parameter $\theta\in\{0,1\}$.
Let $\Gamma^0$, an approximation to $\Gamma(0)$, as well as
$\vec U^{0} \in [S^h_\ell(\Gamma^0)]^3$ and $\vec\kappa^0 \in \qVhz$ be given.
For $m=0,\ldots, M-1$, find $\vec U^{m+1} \in [S^h_\ell(\Gamma^m)]^3$, 
$P^{m+1} \in S^h_{\ell-1}(\Gamma^m)$, 
$\vec X^{m+1}\in\qVh$, $\vec\kappa^{m+1}\in\qVh$ and
$\vec F^{m+1} \in \qVh$ such that
\begin{subequations} \label{eq:qGD}
\begin{align}
&
\rho \left\langle \frac{\vec U^{m+1} - \vec U^m}{\tau}, 
\vec \xi \right\rangle_{\Gamma^m}
+ \tfrac\theta2 \rho \left\langle \nabs\cdot\vec U^m, \vec U^m \cdot \vec\xi 
\right\rangle_{\Gamma^m} 
+ 2\,\mu \left\langle \mat D_s^m (\vec U^{m+1}) , 
\mat D_s^m ( \vec \xi ) \right\rangle_{\Gamma^m}
\nonumber \\ & \qquad\quad
- \left\langle P^{m+1} , 
\nabs\cdot \vec \xi \right\rangle_{\Gamma^m}
= \left\langle \vec g^m, \vec \xi \right\rangle_{\Gamma^m} 
+ \alpha \left\langle\vec F^{m+1}, \vec \xi \right\rangle_{\Gamma^m}
\qquad \forall\ \vec\xi \in \qVh \,, \label{eq:qGDa}\\
& \left\langle \nabs\cdot\vec U^{m+1}, \eta \right\rangle_{\Gamma^m}  = 0 
\qquad \forall\ \eta \in S^h_{\ell-1}(\Gamma^m)\,,
\label{eq:qGDc} \\
& \left\langle \frac{\vec X^{m+1} - \vec\id}{\tau} ,
\vec\chi \right\rangle_{\Gamma^m}
= \left\langle \vec U^{m+1}, \vec\chi \right\rangle_{\Gamma^m}
 \qquad\forall\ \vec\chi \in \qVh\,,
\label{eq:qGDd} \\
& \left\langle \vec\kappa^{m+1} , \vec\eta \right\rangle_{\Gamma^m}
+ \left\langle \nabs\,\vec X^{m+1}, \nabs\,\vec \eta \right\rangle_{\Gamma^m}
 = 0  \qquad\forall\ \vec\eta \in \qVh\,,\label{eq:qGDe} \\
& \left\langle \vec F^{m+1}, \vec\chi \right\rangle_{\Gamma^m} =
\left\langle \nabs\,\vec\kappa^{m+1} , \nabs\,\vec\chi \right\rangle_{\Gamma^m}
+ \left\langle \nabs\cdot(\vec\Pi_{m-1}^{m}\,\vec\kappa^{m}), \nabs\cdot\vec\chi 
\right\rangle_{\Gamma^m}
 \nonumber \\ & \hspace{2cm}
+\tfrac12 \left\langle |\vec\Pi_{m-1}^{m}\,\vec\kappa^m|^2\,\nabs\,\vec \id,
 \nabs\,\vec\chi \right\rangle_{\Gamma^m}
-2 \left\langle [\nabs\,(\vec\Pi_{m-1}^{m}\,\vec\kappa^{m})]^T, 
 \mat D_s^m(\vec\chi)\,(\nabs\,\vec\id)^T \right\rangle_{\Gamma^m}
\nonumber \\ & \hspace{10cm}
\qquad \forall\ \vec\chi \in \qVh\,, \label{eq:qGDf} 
\end{align}
\end{subequations}
and set $\Gamma^{m+1} = \vec X^{m+1}(\Gamma^m)$. 

We observe that \eqref{eq:qGD} is a linear scheme in that
it leads to a linear system of equations for the unknowns 
$(\vec U^{m+1}, P^{m+1}, \vec X^{m+1}, \vec \kappa^{m+1}, \vec F^{m+1})$ 
at each time level. Moreover, in the case $\alpha=0$ the system
(\ref{eq:qGDa}--c) for $(\vec U^{m+1}, P^{m+1}, \vec X^{m+1})$
decouples from (\ref{eq:qGDe},e), and $(\vec \kappa^{m+1}, \vec F^{m+1})$ need
not be computed at all.
Observe that the implicit/explicit treatment on the right hand side of
\eqref{eq:qGDf} follows the one from \cite[(5.3)]{nsns} and allows a
well-posedness proof.

Similarly to \eqref{eq:LBB}, we also say that 
$(\qVh, \Wh)$ 
satisfy a discrete LBB inf-sup condition if there exists a
$C_0 \in \bR_{>0}$, independent of $\{\sigma^m_j\}_{j=1}^{J}$, such that
\begin{equation} \label{eq:qLBB0}
\inf_{\eta \in\Wh} 
\sup_{\vec \xi \in \qVh}
\frac{\left\langle \eta, \nabs\cdot\vec \xi
\right\rangle_{\Gamma^m}}
{\| \eta \|_{0,\Gamma^m}
\,(\|\vec \xi\|_{0,\Gamma^m} + 
 \| \mat{\mathcal{P}}^m\,\vec\xi \|_{1,\Gamma^m,h})} \geq C_0\,,
\end{equation}
where 
$\| \eta \|_{0,\Gamma^m}^2 =
\left\langle \eta,\eta \right\rangle_{\Gamma^m}$
and $\| \vec\eta \|_{1,\Gamma^m,h}^2 =
\left\langle \vec\eta,\vec\eta \right\rangle_{\Gamma^m}
+ \sum_{j = 1}^{J} \int_{\sigma_j^m} |\nabs\,\vec\eta|^2 \dH{2}$.
We note that the condition \eqref{eq:qLBB0} is similar to
an LBB condition for $P\ell$--$P(\ell-1)$ elements in Euclidean space, 
which is known to hold.
Unfortunately, it does not seem to be easily possible to extend the ideas in
the proof of Lemma~\ref{lem:LBB} to the discrete setting. For example,
the projections $\mat{\mathcal{P}}^m$ are discontinuous on $\Gamma^m$.
For completeness we mention that in the stationary case, discrete LBB
conditions have recently been proved in 
\cite{HarderingP25,Reusken25}. 

In the absence of the LBB condition (\ref{eq:qLBB0}) we need to consider the
reduced system (\ref{eq:qGDa},c--e), where $\qVh$ in 
(\ref{eq:qGDa}) is replaced by $\uspace^m_0$. Here we define
\begin{align*} 
\uspace^m_0 = 
\left\{ \vec U \in \qVh : 
\left\langle \nabs\cdot \vec U, \eta 
\right\rangle_{\Gamma^m}  = 0 \ \ \forall\ \eta \in \Wh
 \right\} \,.
\end{align*}

\begin{theorem} \label{thm:qGD}
Let $\rho > 0$, $\mu>0$ and let the LBB condition \eqref{eq:qLBB0} hold. 
Then there exists a unique solution $(\vec U^{m+1}, P^{m+1}, \vec X^{m+1}, 
\vec\kappa^{m+1}, \vec F^{m+1}) 
\in \qVh \times \Wh \times [\qVh]^3$ 
to \eqref{eq:qGD}. 
If \eqref{eq:qLBB0} does not hold, then for $\rho > 0$ there exists a 
unique solution \linebreak
$(\vec U^{m+1}, \vec X^{m+1}, \vec\kappa^{m+1}, \vec F^{m+1}) \in 
\uspace^m_0 \times [\qVh]^3$ to the
reduced system {\rm (\ref{eq:qGDa},c--e)} with 
$\qVh$ replaced by $\uspace^m_0$.
\end{theorem}
\begin{proof}
As the system \eqref{eq:qGD} is linear, existence follows from uniqueness.
In order to establish the latter, we consider the homogeneous system.
Find $(\vec U, P, \vec X, \vec\kappa, \vec F) 
\in \qVh \times \Wh \times [\qVh]^3$ such that
\begin{subequations}
\begin{align}
&
 \tfrac1{\tau}\, \rho \left\langle \vec U , \vec\xi
 \right\rangle_{\Gamma^m}
+ 2\,\mu \left\langle \mat D_s^m (\vec U) , 
\mat D_s^m ( \vec \xi ) \right\rangle_{\Gamma^m}
- \left\langle P, 
\nabs\cdot \vec \xi \right\rangle_{\Gamma^m}
- \alpha \left\langle \vec F, \vec \xi \right\rangle_{\Gamma^m}
= 0 
\nonumber \\ & \hspace{9cm}
 \qquad \forall\ \vec\xi \in \qVh \,, \label{eq:qproofa}\\
& \left\langle \nabs\cdot\vec U, \eta 
\right\rangle_{\Gamma^m}  = 0 
\qquad \forall\ \eta \in \Wh\,,
\label{eq:qproofc} \\
& \tfrac1{\tau} \left\langle \vec X,\vec\chi \right\rangle_{\Gamma^m}
= \left\langle \vec U, \vec\chi \right\rangle_{\Gamma^m}
 \qquad\forall\ \vec\chi \in \qVh\,,
\label{eq:qproofd} \\
& \left\langle \vec\kappa , \vec\eta \right\rangle_{\Gamma^m}
+ \left\langle \nabs\,\vec X, \nabs\,\vec \eta \right\rangle_{\Gamma^m}
 = 0  \qquad\forall\ \vec\eta \in \qVh \,,\label{eq:qproofe}  \\
& \left\langle \vec F, \vec\chi \right\rangle_{\Gamma^m} -
\left\langle \nabs\,\vec\kappa , \nabs\,\vec\chi \right\rangle_{\Gamma^m}
= 0 \qquad \forall\ \vec\chi \in \qVh\,. \label{eq:qprooff} 
\end{align}
\end{subequations}
Choosing $\vec\xi=\vec U$ in (\ref{eq:qproofa}),
$\eta = P$ in (\ref{eq:qproofc}), 
$\vec\chi = \vec F$ in (\ref{eq:qproofd}), 
$\vec\eta=\vec\kappa$ in (\ref{eq:qproofe})
and $\vec\chi = \vec X$ in (\ref{eq:qprooff}) 
yields that
\begin{align}
& \rho \left\langle \vec U , \vec U \right\rangle_{\Gamma^m}
+ 2\,\tau\,\mu \left\langle \mat D_s^m (\vec U) , 
\mat D_s^m (\vec U) \right\rangle_{\Gamma^m}
+ \alpha \left\langle \vec\kappa, \vec\kappa \right\rangle_{\Gamma^m}
=0\,. \label{eq:qproof2GD}
\end{align}
Since $\rho > 0$, it immediately follows from (\ref{eq:qproof2GD})
that $\vec U = \vec 0$, and so \eqref{eq:qproofd}, 
\eqref{eq:qproofe} and \eqref{eq:qprooff} yield $\vec X = \vec 0$,
$\vec \kappa = \vec 0$ and $\vec F = \vec0$, in that order.
Moreover, if (\ref{eq:qLBB0}) holds then
(\ref{eq:qproofa}) with $\vec U = \vec 0$ and $\vec F = \vec0$ implies
that $P = 0$.
This shows existence and uniqueness of 
$(\vec U^{m+1}, P^{m+1}, \vec X^{m+1}, \vec\kappa^{m+1},
\vec F^{m+1}) 
\in \qVh\times \Wh \times [\qVh]^3$ 
to \eqref{eq:qGD}.

The proof for the reduced system is very similar. The homogeneous system to
consider is (\ref{eq:qproofa},c--e) with $\qVh$ 
replaced by $\uspace^m_0$. As before, we infer that (\ref{eq:qproof2GD}) holds,
which yields that $\vec U = \vec 0$, 
and hence $\vec X = \vec 0$, $\vec \kappa = \vec 0$ and $\vec F = \vec0$.
\end{proof}

\begin{remark} \label{rem:rho}
The fact that we are only able to prove existence and uniqueness to 
\eqref{eq:qGD} in the case $\rho>0$ is related to the lack of a Korn
inequality on surfaces for nontangential vector fields. 
In particular, it is not possible to infer from
\eqref{eq:qproof2GD} with $\rho=0$ that $\vec U = \vec 0$. This is related to
the observation made below \eqref{eq:diff2} at the beginning of the paper.
\end{remark}

\begin{remark} \label{rem:discconst}
Let $\vec g^0 = \vec 0$, $\alpha=0$ and 
$\vec U^0 = \vec u_0 = \vec b \in \bR^3$. Then
$(\vec U^1, P^1, \vec X^1) = (\vec b, 0, 
\vec\id_{\mid_{\Gamma^0}} + \tau\,\vec b)$ solves 
\eqref{eq:qGDa}, \eqref{eq:qGDc} 
and \eqref{eq:qGDd}. That means for the trivial problem of a translating
surface, our scheme produces the exact solution, in the sense that the discrete
surface is translated exactly. See also Remark~\ref{rem:judge}.
\end{remark}

\setcounter{equation}{0}
\section{Solution methods} \label{sec:sol}

On introducing the obvious notation, 
the linear system \eqref{eq:qGD} can be written as
\begin{equation}
\begin{pmatrix}
 \mat B & \vec{\mathcal{C}} & 0 & 0 & -\alpha\,\mat M_{F,U} \\
 \vec {\mathcal{C}}^T & 0 & 0 & 0 & 0 \\
 \mat M_{F,U}^T & 0 & 0 & -\frac1{\tau}\,\mat{M} & 0 \\
0 & 0 & \mat{M} & \mat{A} & 0 \\
0 & 0 & -\mat A & 0 & \mat M 
\end{pmatrix} 
\begin{pmatrix} \vec U^{m+1} \\ P^{m+1} \\ \vec\kappa^{m+1} \\ 
\delta\vec{X}^{m+1} \\ \vec F^{m+1} \end{pmatrix}
=
\begin{pmatrix} \vec b \\ 0 \\
0 \\ -\mat{A}\,\vec{X}^{m} \\ 
 \mat{\mathcal{Z}}\,\vec\kappa^{m} + \mat A_{\vec\kappa}\,\vec
X^{m} \end{pmatrix} \,,
\label{eq:lin}
\end{equation}
where we have omitted the dependence of the matrices and the vector
$\vec b$ on $m$ for notational simplicity.
We refer to \cite{nsns,bgnreview} for more details corresponding to closely
related linear systems of equations.
For the solution of (\ref{eq:lin}) a Schur complement approach similar to
\cite{fluidfbp} can be used. In particular, the Schur approach eliminates 
$(\vec\kappa^{m+1},\delta \vec X^{m+1},\vec F^{m+1})$ from 
(\ref{eq:lin}) as follows. Let 
\begin{equation*} 
\Theta= \begin{pmatrix}
 0 & - \frac1{\tau}\,\mat{M} & 0\\
\mat{M} & \mat{A} & 0 \\
-\mat A & 0 & \mat M 
\end{pmatrix} \,.
\end{equation*}
Then (\ref{eq:lin}) can be reduced to
\begin{subequations}
\begin{align} \label{eq:schur}
&
\begin{pmatrix}
\mat B + \alpha\,\mat T
& \vec {\mathcal{C}} \\
\vec {\mathcal{C}}^T & 0 
\end{pmatrix}
\begin{pmatrix}
\vec U^{m+1} \\ P^{m+1} 
\end{pmatrix}
= \begin{pmatrix}
\vec b + \alpha\,\vec c \\
0
\end{pmatrix}
\end{align}
and
\begin{equation} \label{eq:schurb}
\begin{pmatrix}
\vec\kappa^{m+1} \\ \delta\vec{X}^{m+1} \\ \vec F^{m+1}
\end{pmatrix}
 = \Theta^{-1}\,
\begin{pmatrix}
-\mat M_{F,U}^T\,\vec U^{m+1} \\ -\mat{A}\,\vec{X}^{m} \\
\mat{\mathcal{Z}}\,\vec\kappa^{m} + \mat A_{\vec\kappa}\,\vec X^m
\end{pmatrix}
\,.
\end{equation}
\end{subequations}
In (\ref{eq:schur}) we have used the definitions
$$\mat T = (0\ 0\ \mat M_{F,U})\,\Theta^{-1}\,
\begin{pmatrix} \mat M_{F,U}^T \\ 0 \\ 0 \end{pmatrix}
= \tau\,\mat M_{F,U}\,\mat M^{-1}\,\mat A\,\mat M^{-1}\,
\mat A\,\mat M^{-1}\,\mat M_{F,U}^T$$
and 
\begin{align*}
\vec c & = (0\ 0\ \mat M_{F,U})\,\Theta^{-1}\,
\begin{pmatrix}
0 \\ - \mat{A}\,\vec{X}^{m} \\
\mat{\mathcal{Z}}\,\vec\kappa^{m} + \mat A_{\vec\kappa}\,\vec X^m
\end{pmatrix}
\\ &
= \mat M_{F,U}\,\mat M^{-1}\,[\mat{\mathcal{Z}}\,\vec\kappa^m + 
\mat A_{\vec\kappa}\,\vec X^m
- \mat A\,\mat M^{-1}\,\mat{A}\,\vec{X}^m]\,.
\end{align*}
For the linear system (\ref{eq:schur})
well-known solution methods for finite element discretizations for the 
standard Navier--Stokes equations may be employed. We refer to 
\cite[\S5]{fluidfbp}, where we describe such solution methods in 
detail for a similar situation.

\setcounter{equation}{0} 
\section{Numerical results} \label{sec:nr2}
We implemented our fully discrete finite approximation 
\eqref{eq:qGD} for the case $\ell=2$, i.e.\ for piecewise quadratic surfaces,
within the finite element toolbox ALBERTA, see \cite{Alberta}. For all our
computations we will use $\theta=1$.
The inner products in \eqref{eq:qGD} with nonpolynomial integrands are
evaluated with the help of quadrature rules that are exact for polynomials of
degree up to 17.
The linear systems of equations \eqref{eq:lin} arising at each time level we
solve with the Schur complement solver introduced in Section~\ref{sec:sol}.
Here we employ a GMRES iterative solver for \eqref{eq:schur}, with
preconditioner a least square solution for the block matrix on the left hand
side with $\alpha=0$. For this least square solution we employ the sparse
factorization package SPQR, see \cite{Davis11}, while for the computation
of \eqref{eq:schurb} we use the package UMFPACK, see \cite{Davis04}.

Given the initial surface $\Gamma^0$, we define the initial data
$\vec\kappa^0 \in \qVhz$ as the solution to
\begin{equation*} 
 \left\langle \vec\kappa^{0} , \vec\eta \right\rangle_{\Gamma^0}
+ \left\langle \nabs\,\vec\id, \nabs\,\vec \eta \right\rangle_{\Gamma^0}
 = 0   \qquad\forall\ \vec\eta \in \qVhz\,.
\end{equation*}
Throughout this section, unless stated
otherwise, we set $\tau=10^{-3}$, fix 
$\vec u_0 = \vec g = \vec 0$,
$\rho = \mu = \alpha = 1$, and consider the case of the bending energy 
\eqref{eq:bendingE}.
At times we will discuss the discrete energies of the numerical solutions, 
which, on recalling Theorem~\ref{thm:stabqGD} are defined by
$$
E^{m+1} = \tfrac12 \rho \left\langle \vec U^{m+1}, \vec U^{m+1} 
\right\rangle_{\Gamma^m} + \alpha E_\alpha^{m+1},\quad
E_\alpha^{m+1} = \tfrac12\left\langle \vec\kappa^{m+1}, \vec\kappa^{m+1} 
\right\rangle_{\Gamma^m} .
$$

\newcommand{\errorXxT}{\|\vec{X}(T) - \vec{x}(T)\|_{L^\infty}}
\newcommand{\errorGgT}{\|\Gamma^h(T) - \Gamma(T)\|_{L^\infty}}
\newcommand{\errorPpT}{\|P(T) - p(T)\|_{L^2}}
\newcommand{\errorGg}{\|\Gamma^h - \Gamma\|_{L^\infty}}
\newcommand{\errorPp}{\|P - p\|_{L^2(L^2)}}
We begin with a convergence experiment for a family of radially
expanding and shrinking spheres. In particular, we show in
Appendix~\ref{sec:App_truesol} that
a family of spheres $\Gamma(t) = r(t)\,\bS^{2}$, 
with outer unit normal $\vec\nu = \frac{\vec\id}{|\vec\id|}$, 
together with $\vec u(\cdot, t) = r'(t)\,\vec\nu$ and
\begin{equation*} 
p(t) = \tfrac12\rho\,r''(t)\,r(t) + 2\,\mu\,\frac{r'(t)}{r(t)} 
- \tfrac12 \alpha\,f_\Gamma(t)\,r(t)\,,
\end{equation*}
is a solution to \eqref{eq:esns} with the second equation replaced by
$\nabs\cdot\vec u = 2\,\frac{r'(t)}{r(t)}$. We note that in the case of the
bending energy \eqref{eq:bendingE} it holds that
$f_\Gamma(t) = 0$.
We use this solution to perform some convergence tests.
When we choose $\rho=\mu=\alpha=1$ and fix
$r(t) = 1 + \tfrac12\,\sin (2\,\pi\,t)$, we obtain the errors reported in
Table~\ref{tab:q3dsin}. Here we have defined the error quantity
\begin{equation*}
    \errorGg =
    \max_{m=1,\ldots, M} 
    \max_{k=1,\ldots,K}
    \operatorname{dist}(\vec{q}^m_{k}, \Gamma(t_m)) .
\end{equation*}
The numbers in Table~\ref{tab:q3dsin} indicate a convergence order of
$\mathcal{O}(h^3)$. When we repeat this experiment for the value $\alpha=0$,
the discrete surfaces exhibit oscillations and no convergence to the exact 
solution can be observed.
\begin{table}
\center
\begin{tabular}{|r|c|c|c|c|c|}
\hline
$J$  & $h_0$ & $\errorGg$ & EOC & $\errorPpT$ & EOC 
\\ \hline
384  & 4.0994e-01 & 2.3137e-01& ---  & 6.2949e-00& ---  \\ 
768  & 2.7688e-01 & 7.0352e-02& 3.03 & 9.9219e-01& 4.71 \\ 
1536 & 2.0854e-01 & 2.8472e-02& 3.19 & 3.6524e-01& 3.53 \\ 
3072 & 1.3975e-01 & 8.3914e-03& 3.05 & 1.1255e-01& 2.94 \\ 
6144 & 1.0472e-01 & 3.5054e-03& 3.02 & 4.6171e-02& 3.09 \\ 
12288& 7.0041e-02 & 1.0503e-03& 3.00 & 1.4754e-02& 2.84 \\ 
24576& 5.2416e-02 & 4.4005e-04& 3.00 & 6.2182e-03& 2.98 \\ 
\hline
\end{tabular}
\caption{\eqref{eq:qGD}: 
$\rho=\mu=\alpha=1$, $r(t) = 1 + \tfrac12\,\sin (2\,\pi\,t)$. 
$T=1$, $\tau= h_0^{3}$.}
\label{tab:q3dsin}
\end{table}%

\begin{figure}
\center
\mbox{
\includegraphics[angle=-0,width=0.2\textwidth]{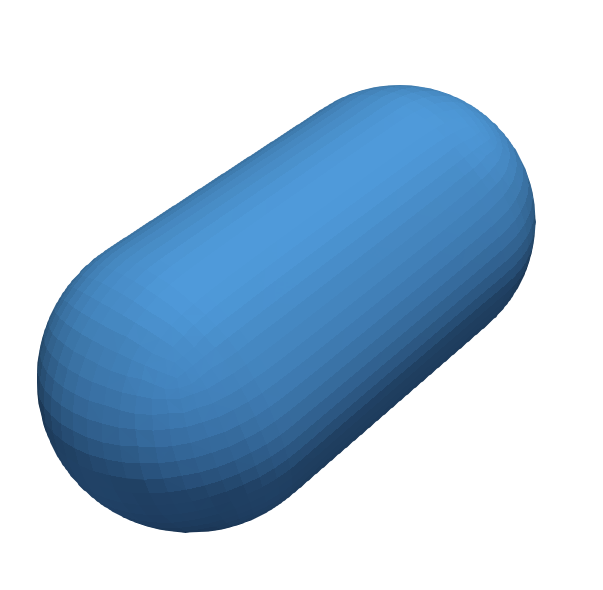}
\includegraphics[angle=-0,width=0.2\textwidth]{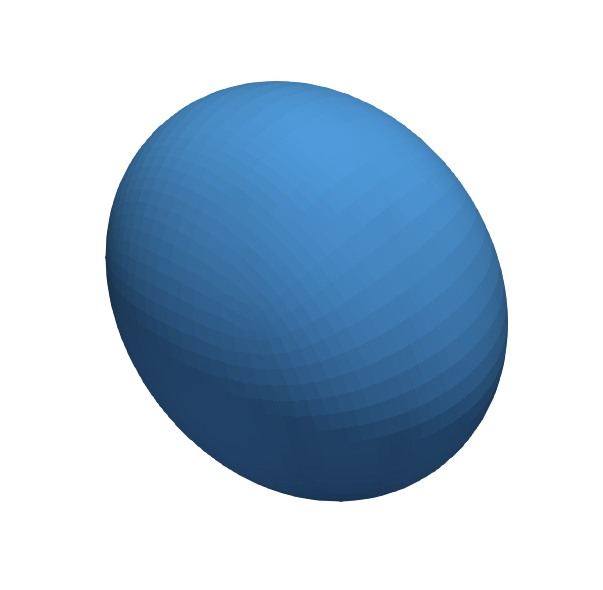}
\includegraphics[angle=-0,width=0.2\textwidth]{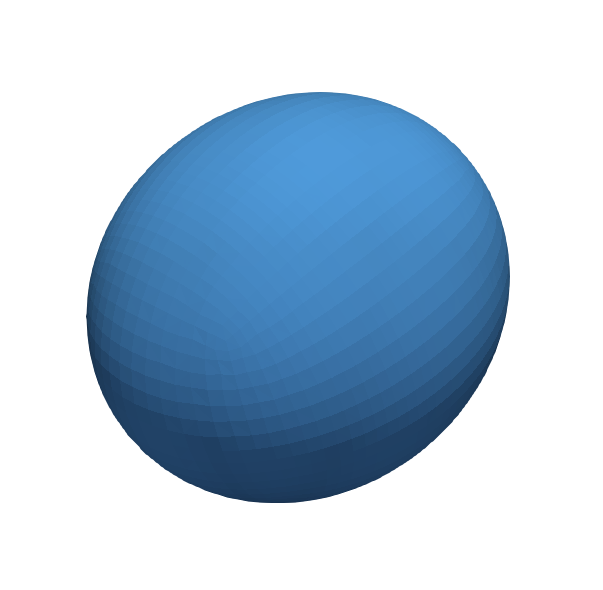}
\includegraphics[angle=-0,width=0.2\textwidth]{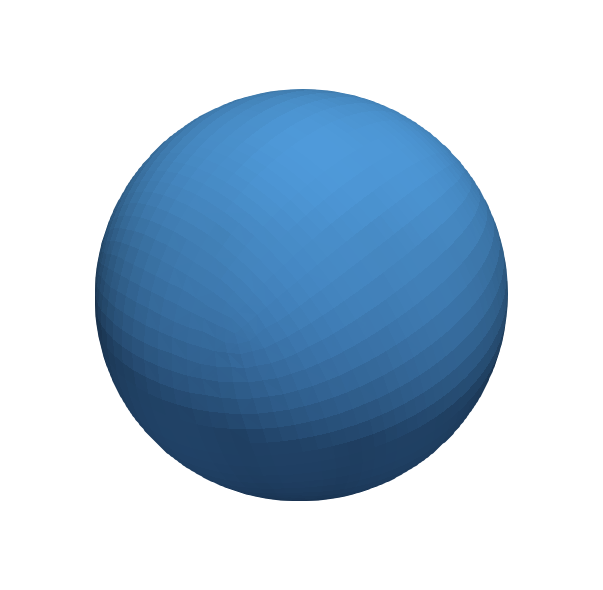}
}
\includegraphics[angle=-90,width=0.45\textwidth]{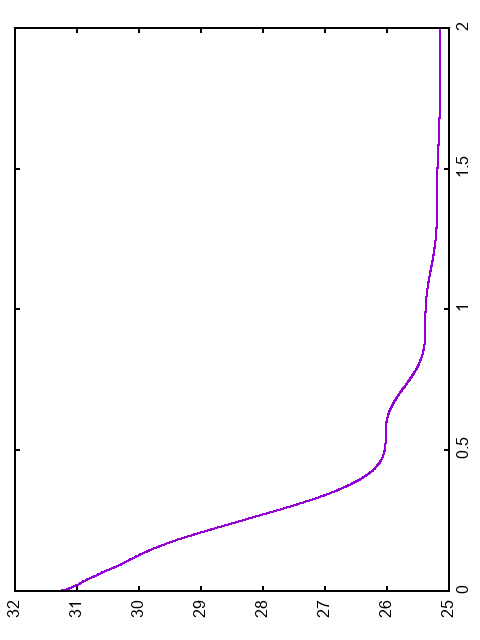}
\includegraphics[angle=-90,width=0.45\textwidth]{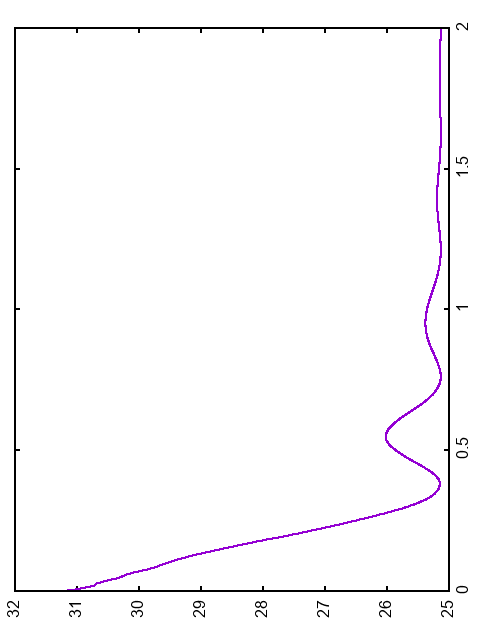}
\caption{($\rho=\mu=\alpha=1$)
The surface $\Gamma^m$ at times $t=0,0.6,1,2$. 
Below we show plots of $E^{m}$ (left) and $E^{m}_\alpha$ (right) over time.}
\label{fig:cig211}
\end{figure}

\begin{figure}
\center
\mbox{
\includegraphics[angle=-0,width=0.2\textwidth]{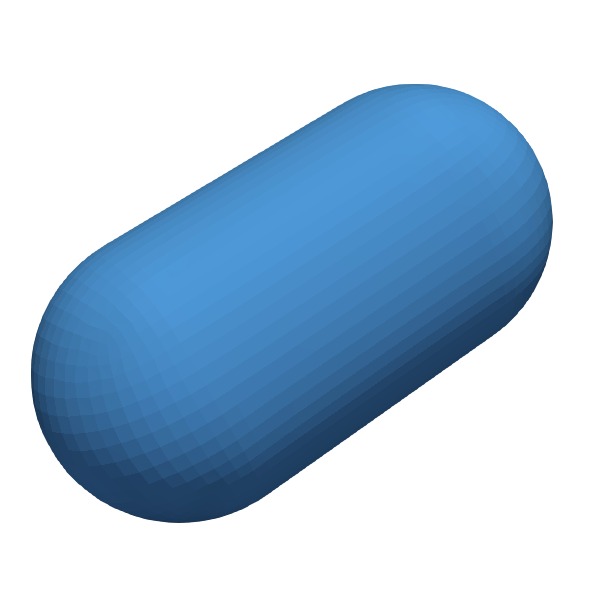}
\includegraphics[angle=-0,width=0.2\textwidth]{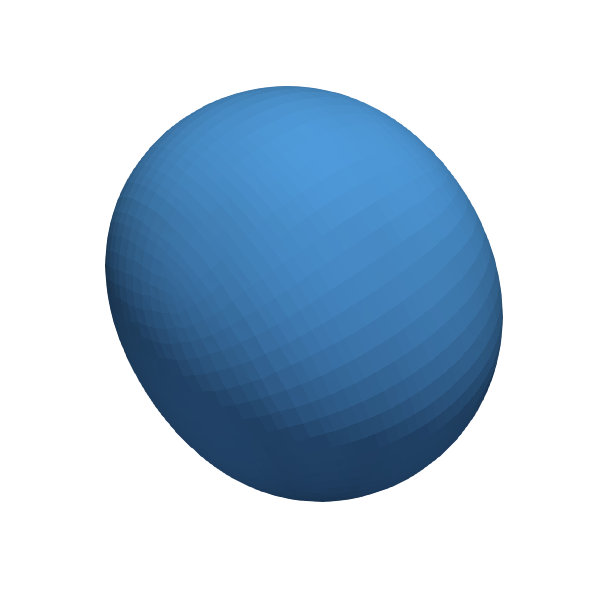}
\includegraphics[angle=-0,width=0.2\textwidth]{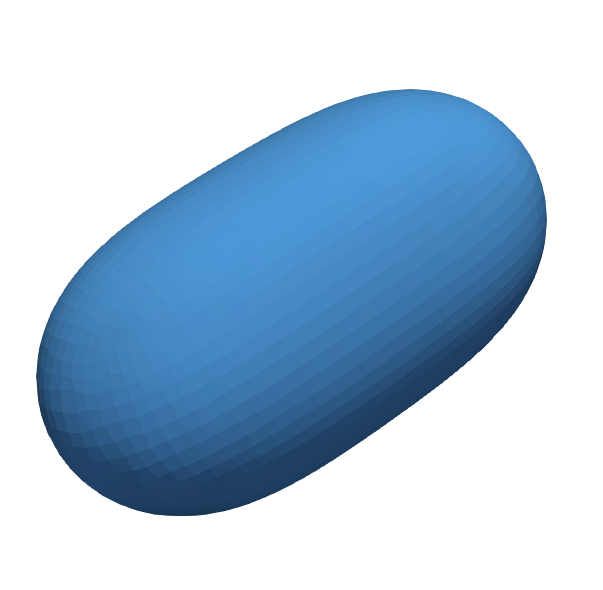}
\includegraphics[angle=-0,width=0.2\textwidth]{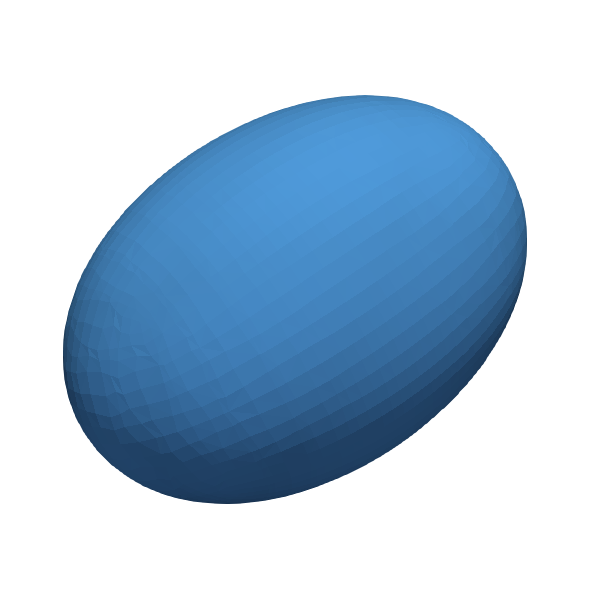}
}
\includegraphics[angle=-90,width=0.45\textwidth]{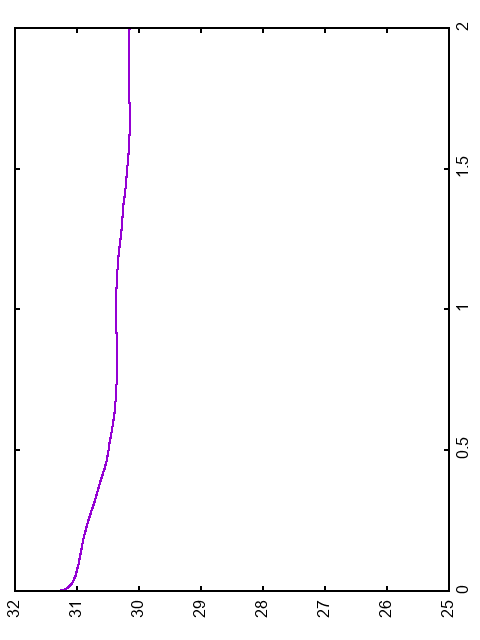}
\includegraphics[angle=-90,width=0.45\textwidth]{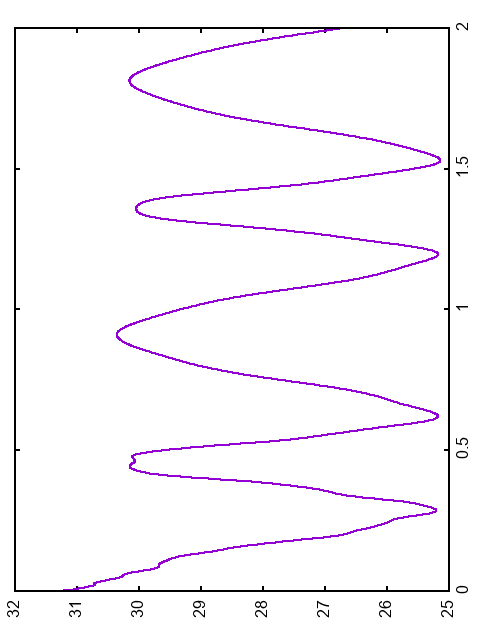}
\caption{($\rho=10^{4}\cdot\mu=\alpha=1$)
The surface $\Gamma^m$ at times $t=0,0.6,1,2$. 
Below we show plots of $E^{m}$ (left) and $E^{m}_\alpha$ (right) over time.}
\label{fig:cig211mu1e-4}
\end{figure}

We now investigate the evolution for an initially stationary 
$2\times1\times1$ tubular shape. For an approximation that uses
$J=1280$ elements, we visualize some snapshots of the evolving surface
together with plots over time of the total discrete energy $E^{m}$
and the discrete bending energy $E^{m}_\alpha$ in
Figure~\ref{fig:cig211}.
When we reduce the interfacial shear viscosity to $\mu=10^{-4}$, then the flow
is far less dissipative and the evolving surface oscillates much more and for
much longer. We demonstrate that with the corresponding plots in
Figure~\ref{fig:cig211mu1e-4}.

\begin{figure}
\center
\mbox{
\includegraphics[angle=-0,width=0.4\textwidth]{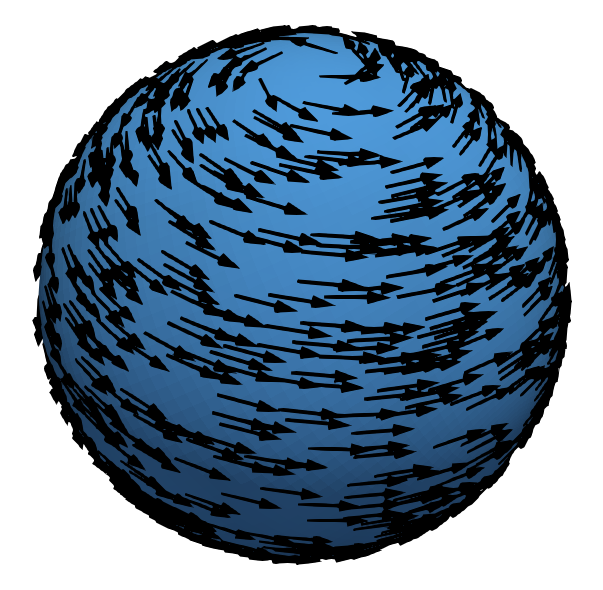}
\includegraphics[angle=-0,width=0.4\textwidth]{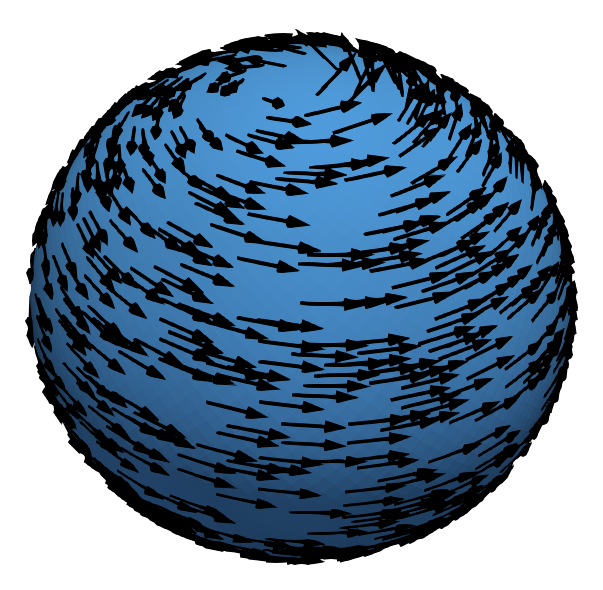}
}
\includegraphics[angle=-90,width=0.45\textwidth]{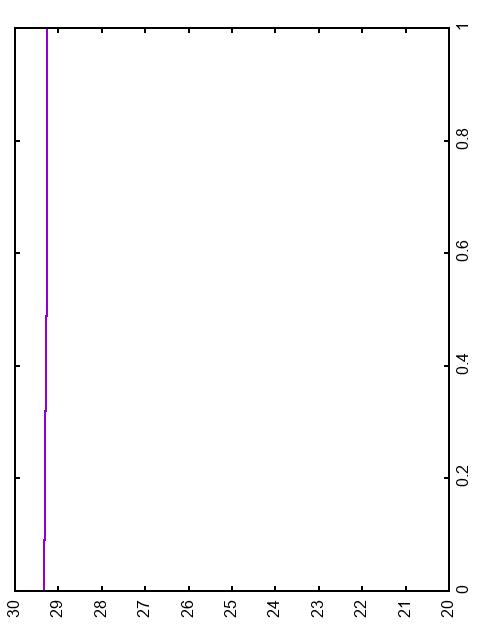}
\includegraphics[angle=-90,width=0.45\textwidth]{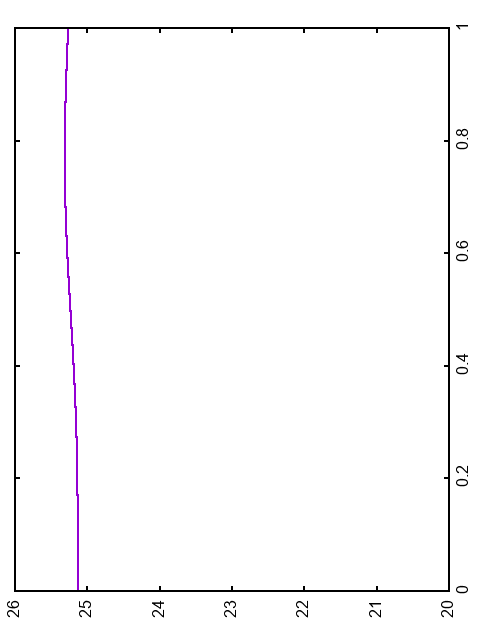}
\caption{($\rho=\mu=\alpha=1$)
The surface $\Gamma^m$, together with a visualization of the velocity field
$\vec U^m$ at times $t=0,1$. 
Below we show plots of $E^{m}$ (left) and $E^{m}_\alpha$ (right) over time.}
\label{fig:killingfieldalpha1}
\end{figure}

We also consider an experiment for a Killing field, recall
Remark~\ref{rem:killingfields}. In particular, we consider the unit sphere
$\bS^2$ with the initial velocity $\vec u_0(\vec z) = (-z_2, z_1, 0)^T$
for $\vec z \in \bS^2$, so that the solution to \eqref{eq:esns} with
\eqref{eq:fGamma} is given by
$(\Gamma(t), \vec u(\cdot, t), p(t)) = (\bS^2, \vec u_0, 0)$.
For the discretization parameters we use $J=4096$ and $\tau = 10^{-4}$.
It can be seen in Figure~\ref{fig:killingfieldalpha1} that the sphere rotates
as expected. Moreover, the fully discrete finite element scheme is not
able to maintain a perfectly spherical approximation of the rotating surface,
which becomes slightly flatter as time goes on. This behaviour can also be
noticed in the two energy plots, where we can see a slightly oscillatory
behaviour of the bending energy plot. This indicates that the discrete surface
undergoes flattening and expanding as time goes on.

\begin{figure}
\center
\includegraphics[angle=-0,width=0.3\textwidth]{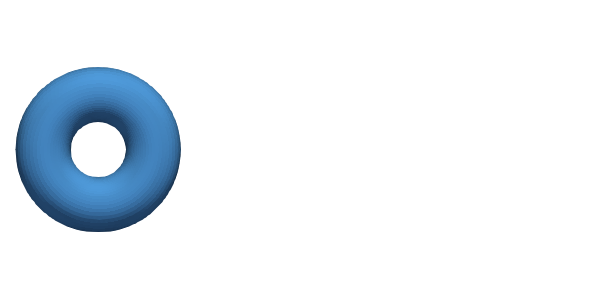}
\includegraphics[angle=-0,width=0.3\textwidth]{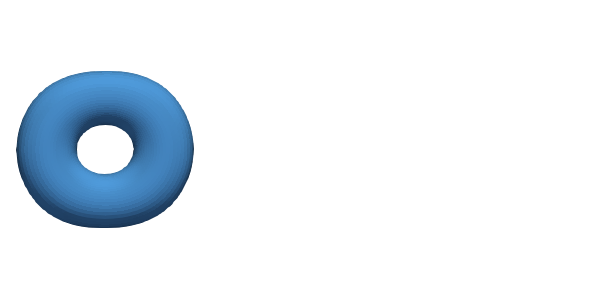}
\includegraphics[angle=-0,width=0.3\textwidth]{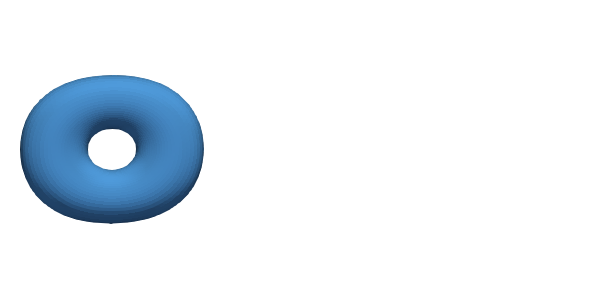}
\includegraphics[angle=-0,width=0.3\textwidth]{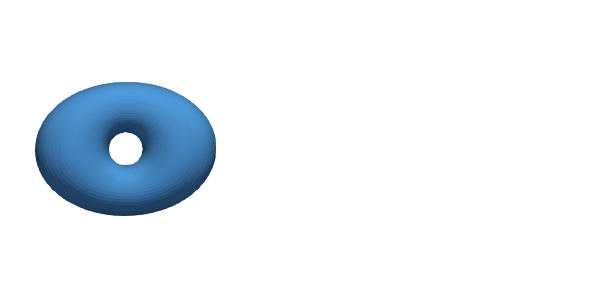}
\includegraphics[angle=-0,width=0.3\textwidth]{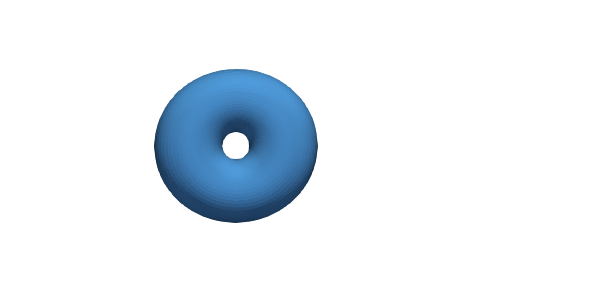}
\includegraphics[angle=-0,width=0.3\textwidth]{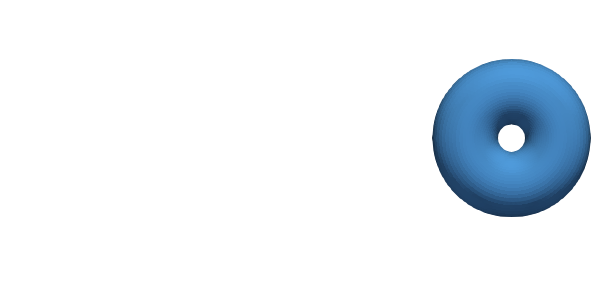}
\includegraphics[angle=-90,width=0.45\textwidth]{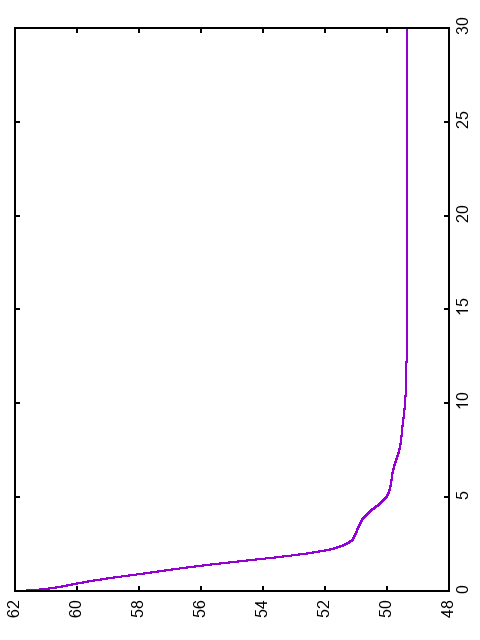}
\includegraphics[angle=-90,width=0.45\textwidth]{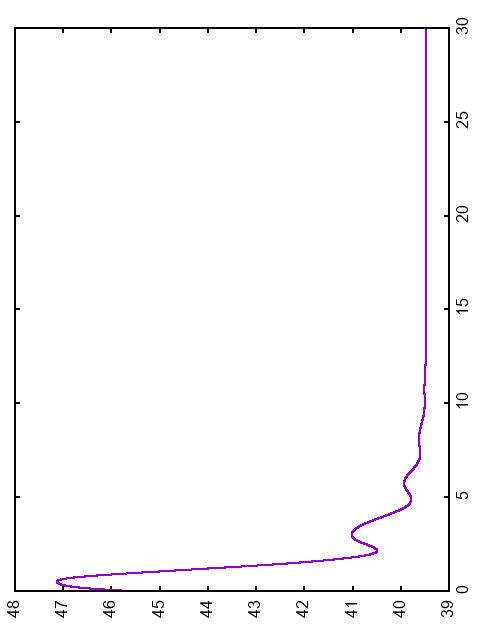}
\caption{($\rho=\mu=\alpha=1$)
The surface $\Gamma^m$ at times $t=0,0.5,1,2,10,30$. 
Below we show plots of $E^{m}$ (left) and $E^{m}_\alpha$ (right) over time.}
\label{fig:torusalpha1}
\end{figure}

We end this section with a simulation for a torus that moves due to a given 
nonzero initial velocity $\vec u_0$. For the initial torus $\Gamma_0$,
centred at the origin, we choose
the radii $R=2$ and $r=1$, while for the initial velocity we consider
$\vec u_0(\vec z) = \begin{cases} 1 & z_1 \geq 0 \\ 0 & z_1 < 0 \end{cases}$.
As a consequence, the moving surface initially moves and stretches in the
$x$-direction, and eventually it seems to settle on a static shape, which
appears to be close to the Clifford torus, that travels with a constant
velocity in the $x$-direction. See Figure~\ref{fig:torusalpha1} for a
visualization of the numerical results, where we used $J=2760$ elements
to represent the discrete surfaces.

\begin{figure}
\center
\mbox{
\includegraphics[angle=-0,width=0.25\textwidth]{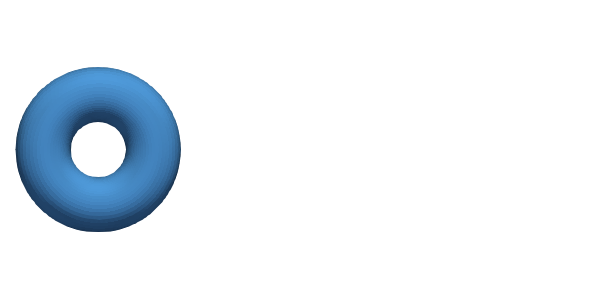}
\includegraphics[angle=-0,width=0.25\textwidth]{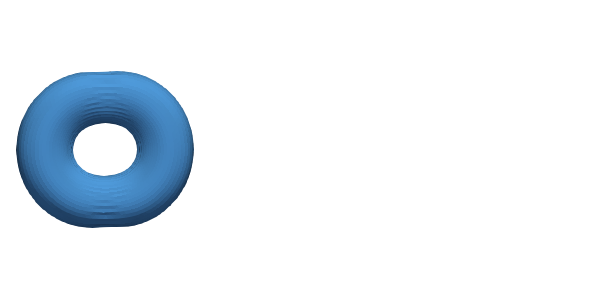}
\includegraphics[angle=-0,width=0.25\textwidth]{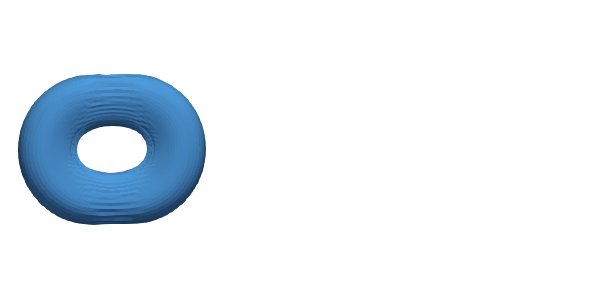}
\includegraphics[angle=-0,width=0.25\textwidth]{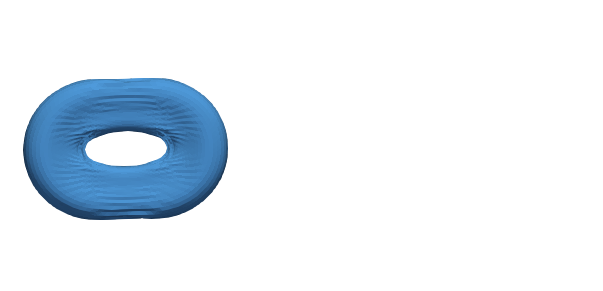}
}
\includegraphics[angle=-90,width=0.45\textwidth]{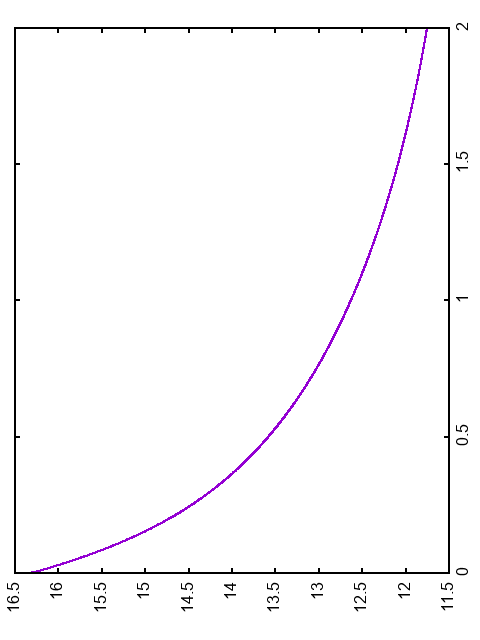}
\includegraphics[angle=-90,width=0.45\textwidth]{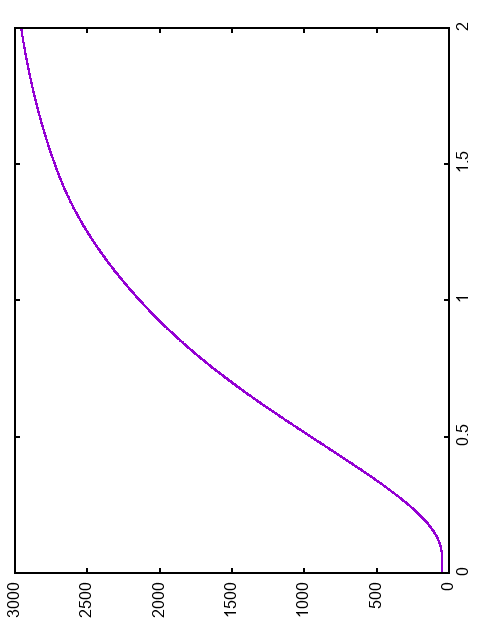}
\caption{($\rho=\mu=1$, $\alpha=0$)
The surface $\Gamma^m$ at times $t=0,0.5,1,2$. 
Below we show plots of $E^{m}$ (left) and $E^{m}_\alpha$ (right) over time.}
\label{fig:torusalpha0}
\end{figure}

We are also interested in the special case $\alpha=0$. In all our numerical
experiments with $\alpha=0$
we experienced strong instabilities and oscillations in the discrete surfaces.
By way of example, we show some snapshots of a repetition of the experiment in
Figure~\ref{fig:torusalpha0}, but now for $\alpha=0$. At the final time we can
see the onset of some oscillations, which would get worse if the 
simulation was continued.

\begin{figure}
\center
\mbox{
\includegraphics[angle=-0,width=0.25\textwidth]{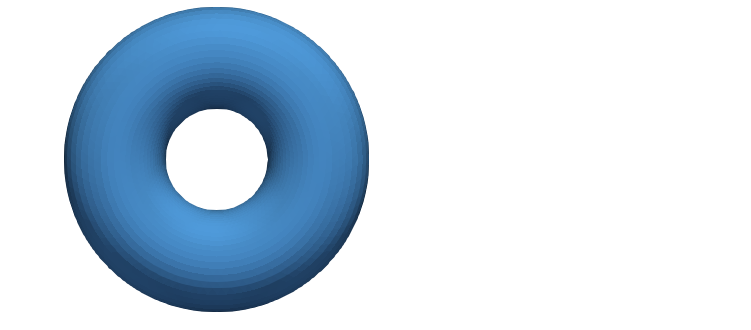}
\includegraphics[angle=-0,width=0.25\textwidth]{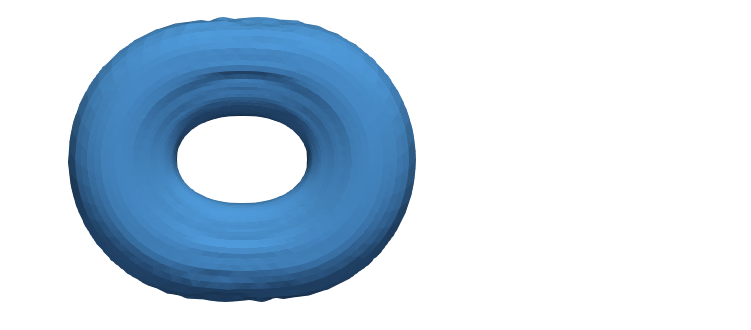}
\includegraphics[angle=-0,width=0.25\textwidth]{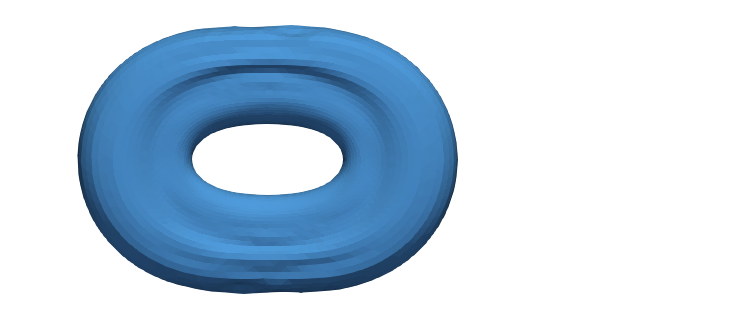}
\includegraphics[angle=-0,width=0.25\textwidth]{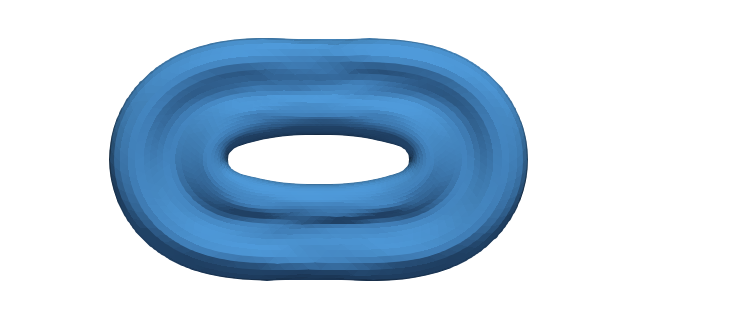}
}
\includegraphics[angle=-90,width=0.45\textwidth]{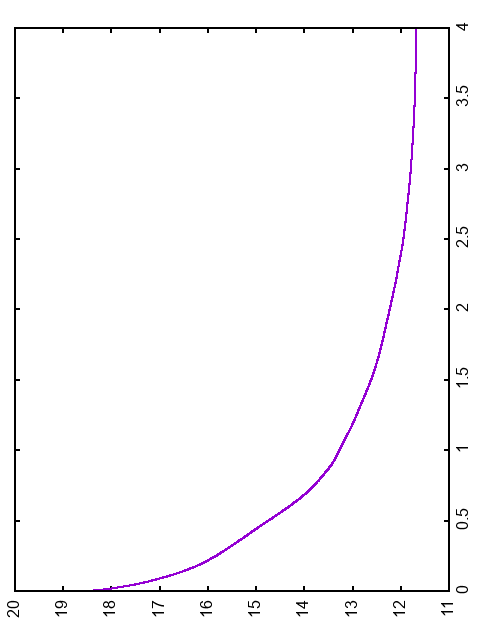}
\includegraphics[angle=-90,width=0.45\textwidth]{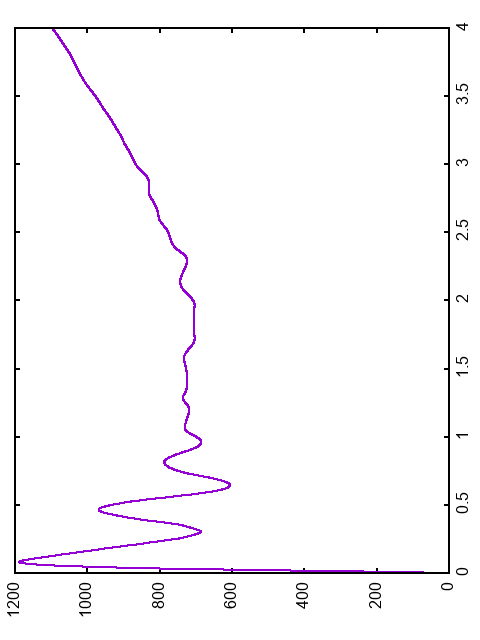}
\caption{($\rho=\mu=1$, $\alpha=10^{-3}$)
The surface $\Gamma^m$ at times $t=0,1,2,4$. 
Below we show plots of $E^{m}$ (left) and $E^{m}_\alpha$ (right) over time.}
\label{fig:torusalpha1e-3}
\end{figure}

If we choose $\alpha$ small but positive, on the other hand, then the evolution
appears to be smooth for much longer. 
A simulation for $\alpha=10^{-3}$ is shown in Figure~\ref{fig:torusalpha1e-3}.
Soon after strong oscillations appear, which we believe are caused by numerical
issues, rather than by the underlying flow.

\newpage
\begin{appendix}
\renewcommand{\theequation}{\Alph{section}.\arabic{equation}}

\setcounter{equation}{0} 
\section{Separation of normal and tangential part} \label{sec:App_JOR}

\begin{lemma} \label{lem:JOR18}
It holds for any scalar field $w$ that
\begin{subequations} \label{eq:JOR18}
\begin{equation} \label{eq:JOR18a}
\nabs\cdot(w\,\nabs\,\vec\nu) = - w\,\nabs\,\varkappa 
- w\,|\nabs\,\vec\nu|^2\,\vec\nu + (\nabs\,\vec\nu)\,\nabs\,w\,.
\end{equation}
Moreover, for a tangential vector field $\vec u = \mat{\mathcal{P}}\,\vec u$
it holds that
\begin{align} \label{eq:JOR18b}
\nabs\cdot(\mat{\mathcal{P}}\, \nabs\,\vec u) & =
\Delta_\subH\,\vec u  - (\nabs\,\vec u : \nabs\,\vec\nu)\,\vec\nu\,,\\
\nabs\cdot[(\nabs\,\vec u)^T\, \mat{\mathcal{P}}] & = 
-\varkappa\,(\nabs\,\vec\nu)\,\vec u
- (\nabs\,\vec\nu)^2\,\vec u
- (\nabs\,\vec u : \nabs\,\vec\nu)\,\vec\nu
+ \nabs\,(\nabs\cdot\vec u)\,, \label{eq:JOR18c}
\end{align}
\end{subequations}
where $\Delta_\subH\,\vec u = 
\mat{\mathcal{P}}\,[\nabs\cdot(\mat{\mathcal{P}} \nabs\,\vec u)]$ denotes the
Bochner Laplacian (also called connection Laplacian) for the tangential vector field $\vec u$.

Finally, writing a possibly nontangential vector field
$\vec u = u_\subN\,\vec\nu + \vec u_\subT$, with $\vec u_\subT =
\mat{\mathcal{P}}\,\vec u$, it holds that
\begin{equation} \label{eq:JOR18d}
\vec\nu\cdot \matpartu\,\vec u = \matpartu\,u_\subN 
- \vec u_\subT \cdot \matpartu\,\vec\nu
\quad{and}\quad
\mat{\mathcal{P}}\, \matpartu\,\vec u = 
u_\subN\,\matpartu\,\vec\nu + \mat{\mathcal{P}}\,\matpartu\,\vec u_\subT
\,.
\end{equation}
\end{lemma}
\begin{proof}
Throughout this proof we use the same notation as in the proof of 
\cite[Lemma~15]{bgnreview}. For example, we write
$\nabs = (\partial_{s_1}, \partial_{s_2}, \partial_{s_3})^T$, 
and $\vec\nu = (\nu_1,\nu_2,\nu_3)^T$. Then we compute
\[
[\nabs\cdot(w\,\nabs\,\vec\nu)]_j
= (\partial_{s_i}\,(w\,\partial_{s_i}\,\nu_j))
= w\,\Delta_s\,\nu_j + \nabs\,w \cdot \nabs\,\nu_j\,,
\]
and so $\nabs\cdot(w\,\nabs\,\vec\nu) = w\,\Delta_s\,\vec\nu + 
(\nabs\,\vec\nu)\,\nabs\,w$. Combining this with 
$\Delta_s\,\vec\nu = - \nabs\,\varkappa - |\nabs\,\vec\nu|^2\,\vec\nu$,
see \cite[Lemma~16]{bgnreview},
gives the desired result \eqref{eq:JOR18a}.

{From} now on let $\vec u = \mat{\mathcal{P}}\,\vec u$ be a 
tangential vector field. Then it holds that
\begin{equation} \label{eq:nPu}
\nabs\cdot(\mat{\mathcal{P}} \nabs\,\vec u)  = 
\mat{\mathcal{P}}\,[\nabs\cdot(\mat{\mathcal{P}} \nabs\,\vec u)]
+ (\vec\nu\cdot [\nabs\cdot(\mat{\mathcal{P}} \nabs\,\vec u)])\,\vec\nu
= \Delta_\subH\,\vec u 
+ (\vec\nu\cdot [\nabs\cdot(\mat{\mathcal{P}} \nabs\,\vec u)])\,\vec\nu\,.
\end{equation}
We compute,
on noting $0 = [\mat{\mathcal{P}}\,\vec\nu]_k = 
\nu_i\,\mathcal{P}_{ik}$, 
that
\begin{align}
\vec\nu\cdot [\nabs\cdot(\mat{\mathcal{P}}\, \nabs\,\vec u)]
& = \nu_i\,\partial_{s_j}\,( \mathcal{P}_{ik}\,\partial_{s_j}\,u_k)
= 
- \partial_{s_j}\,\nu_i\,\mathcal{P}_{ik}\,\partial_{s_j}\,u_k
= - \nabs\,\vec\nu : (\mat{\mathcal{P}}\, \nabs\,\vec u) \nonumber \\ &
= - \nabs\,\vec u : (\mat{\mathcal{P}}\, \nabs\,\vec\nu)  
= - \nabs\,\vec u : \nabs\,\vec\nu\,,
\label{eq:nunPu}
\end{align}
where in the last step we have used
$\mat{\mathcal{P}}\, \nabs\,\vec\nu = 
\mat{\mathcal{P}}\, (\nabs\,\vec\nu)^T = (\nabs\,\vec\nu)^T
=\nabs\,\vec\nu$, recall Lemma~12(ii) and Remark~6(v) from \cite{bgnreview}. 
Combining \eqref{eq:nunPu} and \eqref{eq:nPu} gives the desired result
\eqref{eq:JOR18b}. 

Next we show \eqref{eq:JOR18c}. On noting that
$\vec u \cdot \vec \nu = u_k\,\nu_k = 0$, 
$\vec\nu\cdot\nabs\,f = \nu_j\,\partial_{s_j}\,f = 0$
and $\nabs\cdot\vec\nu = \partial_{s_j}\,\nu_j = -\varkappa$,
we compute
\begin{align*} 
& \partial_{s_j} \left( \partial_{s_i}\,u_k\,(\delta_{kj} - \nu_k\,\nu_j)
\right) \nonumber \\ & \quad
= - (\partial_{s_i}\,u_k)\,\partial_{s_j}\,(\nu_k\,\nu_j)
+ (\partial_{s_j}\,\partial_{s_i}\,u_k)\,(\delta_{kj} - \nu_k\,\nu_j) 
\nonumber \\ & \quad
= - (\partial_{s_i}\,u_k)\left[\nu_j\,\partial_{s_j}\,\nu_k
+ \nu_k\,\partial_{s_j}\,\nu_j \right]
+ \partial_{s_j}\,\partial_{s_i}\,u_j 
- \nu_k\,\nu_j\,\partial_{s_j}\,\partial_{s_i}\,u_k
\nonumber \\ & \quad
= - (\partial_{s_i}\,u_k)\,\nu_k\,\partial_{s_j}\,\nu_j
+ \partial_{s_j}\,\partial_{s_i}\,u_j 
= (\partial_{s_i}\,u_k)\,\nu_k\,\varkappa
+ \partial_{s_j}\,\partial_{s_i}\,u_j \nonumber \\ & \quad
= - \varkappa\,u_k\,\partial_{s_i}\,\nu_k
+ \partial_{s_j}\,\partial_{s_i}\,u_j
= - \varkappa\,[(\nabs\,\vec\nu)\,\vec u]_i
+ \partial_{s_j}\,\partial_{s_i}\,u_j\,.
\end{align*}
In addition, it follows from \cite[Lemma~15]{bgnreview},
on noting $0 = \nabs\,(u_j\,\nu_j) = \nu_j\,\nabs\,u_j +
u_j\,\nabs\,\nu_j$ and on recalling from \cite[Lemma~12]{bgnreview} that
$(\nabs\,\vec\nu)^T = \nabs\,\vec\nu$, that
\begin{align*}
\partial_{s_j}\,\partial_{s_i}\,u_j & 
= \partial_{s_i}\,\partial_{s_j}\,u_j 
+ [(\nabs\,\vec\nu)\, \nabs\,u_j]_i\,\nu_j
- [(\nabs\,\vec\nu)\, \nabs\,u_j)]_j\,\nu_i 
 \nonumber \\ & 
= \partial_{s_i}\,\partial_{s_j}\,u_j 
+ (\nabs\,u_j\cdot \partial_{s_i}\,\vec\nu)\,\nu_j
- (\nabs\,u_j\cdot \partial_{s_j}\,\vec\nu)\,\nu_i 
 \nonumber \\ & 
= \partial_{s_i}\,\partial_{s_j}\,u_j 
+ (\nabs\,\nu_j\cdot \partial_{s_i}\,\vec\nu)\,u_j
- (\nabs\,u_j\cdot \partial_{s_j}\,\vec\nu)\,\nu_i 
 \nonumber \\ & 
= \partial_{s_i}\,\partial_{s_j}\,u_j 
- (\partial_{s_k}\,\nu_j)\,(\partial_{s_i}\,\nu_k)\, u_j
- (\partial_{s_k}\,u_j)\,(\partial_{s_j}\,\nu_k)\,\nu_i
\nonumber \\ & 
= [\nabs\,(\nabs\cdot\vec u)]_i
- [ (\nabs\,\vec\nu)^2\,\vec u]_i
- (\nabs\,\vec u : \nabs\,\vec\nu)\,\nu_i\,.
\end{align*}

Finally, an application of the product rule yields
\begin{equation} \label{eq:JOR18daux}
\matpartu\,\vec u = (\matpartu\,u_\subN)\,\vec\nu 
+ u_\subN\,\matpartu\,\vec\nu + \matpartu\,\vec u_\subT \,.
\end{equation}
On noting that $\vec\nu\cdot \matpartu\,\vec\nu = \frac12\,
\matpartu\,|\vec\nu|^2 = 0$
and $\matpartu\,(\vec u_\subT\cdot\vec\nu) = 0$, 
the desired result \eqref{eq:JOR18d} now follows from \eqref{eq:JOR18daux}.
\end{proof}

\begin{remark} \label{rem:d2}
Let $\vec u = \mat{\mathcal{P}}\,\vec u$ be a tangential vector field.
Then it holds that
\begin{equation} \label{eq:3dJOR18c}
\nabs\cdot[(\nabs\,\vec u)^T\, \mat{\mathcal{P}}] = 
\Gauss\,\vec u - (\nabs\,\vec u : \nabs\,\vec\nu)\,\vec\nu
+ \nabs\,(\nabs\cdot\vec u)\,.
\end{equation}
To see this, similarly to the proof of Lemma~12 in \cite{bgnreview},
we choose an orthonormal basis $\{\onbtau_1,\onbtau_2,\vec\nu\}$ of
$\bR^3$, consisting of eigenvectors of $\nabs\,\vec\nu$ 
with corresponding eigenvalues $-\varkappa_1,-\varkappa_2,0$,
where $\varkappa_1$ and $\varkappa_2$ are the principal curvatures.
Hence 
$(\varkappa\,\nabs\,\vec\nu + (\nabs\,\vec\nu)^2)\,\onbtau_i 
= (- (\varkappa_1 + \varkappa_2)\,\varkappa_i + \varkappa_i^2)\,\onbtau_i
= - \varkappa_1\,\varkappa_2\,\onbtau_i$, and so
\begin{equation} \label{eq:WM}
\varkappa\,\nabs\,\vec\nu + (\nabs\,\vec\nu)^2
= -\Gauss\,\mat{\mathcal{P}}\,,
\end{equation}
with $\Gauss=\varkappa_1\,\varkappa_2$ 
denoting the Gauss curvature. Now \eqref{eq:3dJOR18c} directly follows
from \eqref{eq:JOR18c}, since $\vec u$ is tangential.
\end{remark}

We now want to relate \eqref{eq:esns} to \cite[(3.15)]{JankuhnOR18}. 
It follows from \eqref{eq:nabssigma} that \eqref{eq:esns} is equivalent to
$\mathcal V = \vec u \cdot \vec\nu$ and
\begin{equation} \label{eq:app:esns1}
\rho_\Gamma\,\matpartu\,\vec u 
- 2\,\mu_\Gamma\,\nabs\cdot\mat D_s (\vec u) 
+ \nabs\,p_\Gamma + \varkappa\,p_\Gamma\,\vec\nu
 = \vec g + \alpha\, f_\Gamma\,\vec \nu,\quad
\nabs\cdot\vec u = 0\,.
\end{equation}
Rewriting $\vec u = u_\subN\,\vec\nu + \vec u_\subT$, with $\vec u_\subT =
\mat{\mathcal{P}}\,\vec u$, it follows from \eqref{eq:diff2} that 
\eqref{eq:app:esns1} is equivalent to
\begin{subequations} \label{eq:app:esns2}
\begin{align} \label{eq:app:esns2a}
& \rho\,\matpartu\,\vec u 
- 2\,\mu\,\nabs\cdot
\left[u_\subN\, \nabs\,\vec\nu
+ \tfrac12\,(\mat{\mathcal{P}}\,\nabs\,\vec u_\subT + 
(\nabs\,\vec u_\subT)^T\mat{\mathcal{P}})\right]
+ \nabs\,p + \varkappa\,p\,\vec\nu
 = \vec g+ \alpha\, f_\Gamma\,\vec \nu\,,\\
& \nabs\cdot\vec u_\subT = u_\subN\,\varkappa\,. \label{eq:nabsuT}
\end{align}
\end{subequations}
Combining \eqref{eq:app:esns2a} with \eqref{eq:JOR18} 
then gives 
\begin{align} \label{eq:JORboth}
&
\rho\,\matpartu\,\vec u 
+ 2\,\mu\left( u_\subN\,\nabs\,\varkappa 
+ u_\subN\,|\nabs\,\vec\nu|^2\,\vec\nu - (\nabs\,\vec\nu)\,\nabs\,u_\subN
\right) \nonumber \\ &
- \mu\left( \Delta_\subH\,\vec u_\subT 
- 2\,(\nabs\,\vec u_\subT : \nabs\,\vec\nu)\,\vec\nu
-\varkappa\,(\nabs\,\vec\nu)\,\vec u_\subT
- (\nabs\,\vec\nu)^2\,\vec u_\subT
+ \nabs\,(\nabs\cdot\vec u_\subT)\right) 
+ \nabs\,p + \varkappa\,p\,\vec\nu \nonumber \\ & \quad
= \vec g + \alpha\, f_\Gamma\,\vec\nu \,.
\end{align}
Splitting \eqref{eq:JORboth} into normal and tangential parts, and
combining with \eqref{eq:JOR18d} and \eqref{eq:nabsuT}, yields 
\begin{subequations} \label{eq:JOR18final}
\begin{align}
\rho\,\mat{\mathcal{P}}\,\matpartu\,\vec u_\subT & 
= \mu\left( \Delta_\subH\,\vec u_\subT 
+ \nabs\,(\nabs\cdot\vec u_\subT)
- \varkappa\,(\nabs\,\vec\nu)\,\vec u_\subT
- (\nabs\,\vec\nu)^2\,\vec u_\subT
+ 2\, (\nabs\,\vec\nu)\,\nabs\,u_\subN
\right. \nonumber \\ & \qquad\quad \left.
- 2\,u_\subN\,\nabs\,\varkappa  \right) 
- \rho\,u_\subN\,\matpartu\,\vec\nu - \nabs\,p + \vec g_\subT \,, 
\label{eq:JOR18finala}\\ 
\rho\,\matpartu\,u_\subN & = - 2\,\mu\left( u_\subN\,|\nabs\,\vec\nu|^2
+ \nabs\,\vec u_\subT : \nabs\,\vec\nu \right) 
+ \rho\,u_\subT \cdot \matpartu\,\vec\nu - \varkappa\,p
+ g_\subN + \alpha\,f_\Gamma\,, \label{eq:JOR18finalb} \\ 
\nabs\cdot\vec u_\subT & = u_\subN\,\varkappa\,, \label{eq:JOR18finalc}
\end{align}
\end{subequations}
where we have written $\vec g = g_\subN\,\vec\nu + \vec g_\subT$ with 
$\vec g_\subT = \mat{\mathcal{P}} \, \vec g$.

\begin{remark} \label{rem:JOR18}
We note that, similarly to \cite[Lemma~37]{bgnreview}, and on recalling
$\vec u_\subT\cdot\vec\nu=0$ and $(\nabs\,\vec\nu)^T=\nabs\,\vec\nu$, 
it is possible to write 
\[
\matpartu\,\vec\nu = - (\nabs\,\vec u)^T \,\vec\nu
= - (\vec\nu\otimes\nabs\,u_\subN + u_\subN\,\nabs\,\vec\nu 
+ \nabs\,\vec u_\subT )^T\,\vec\nu
= - \nabs\,u_\subN + (\nabs\,\vec\nu)\,\vec u_\subT\,,
\] 
allowing to replace the two material time derivatives on the right hand sides 
of \eqref{eq:JOR18final} with terms involving only spatial derivatives.
This allows an interpretation of \eqref{eq:esns} as a coupled system of PDEs,
where \eqref{eq:JOR18finalb} describes the evolution of $\mathcal V = u_\subN$,
and hence of $\Gamma(t)$, while \eqref{eq:JOR18finala} describes the evolution
of $\vec u_\subT$.
Following \cite{JankuhnOR18}, it is also possible to utilize 
\eqref{eq:JOR18finalc} 
in order to eliminate derivatives of curvature from \eqref{eq:JOR18final},
e.g.\ the term $\nabs\,\varkappa$ in \eqref{eq:JOR18finalc}.
In particular, it holds that
\begin{align*}
\nabs\,(\nabs\cdot \vec u_\subT) - 2\, u_\subN\,\nabs\,\varkappa &
= 2\left(\nabs(u_\subN\,\varkappa) - u_\subN\,\nabs\,\varkappa\right)
- \nabs\,(\nabs\cdot \vec u_\subT) \\ &
= 2\,\varkappa\,\nabs\,u_\subN - \nabs\,(\nabs\cdot \vec u_\subT)\,,
\end{align*}
and so we can rewrite \eqref{eq:JOR18finala} also as
\begin{align*}
\rho\,\mat{\mathcal{P}}\,\matpartu\,\vec u_\subT & 
= \mu\left( \Delta_\subH\,\vec u_\subT - \nabs\,(\nabs\cdot\vec u_\subT) 
- \varkappa\,(\nabs\,\vec\nu)\,\vec u_\subT
- (\nabs\,\vec\nu)^2\,\vec u_\subT
+ 2\,\varkappa\,\nabs\,u_\subN
\right. \nonumber \\ & \qquad\quad \left.
+ 2\, (\nabs\,\vec\nu)\,\nabs\,u_\subN \right) 
- \rho\,u_\subN\,\matpartu\,\vec\nu 
- \nabs\,p + \vec g_\subT\,,
\end{align*}
which on recalling \eqref{eq:WM} is equivalent to
\begin{align} \label{eq:3dJOR18finala}
\rho\,\mat{\mathcal{P}}\,\matpartu\,\vec u_\subT & 
= \mu\left( \Delta_\subH\,\vec u_\subT
- \nabs\,(\nabs\cdot\vec u_\subT)
+ \Gauss\,\vec u_\subT
+ (2\,\varkappa\,\mat{\mathcal{P}} + \nabs\,\vec\nu)\,\nabs\,u_\subN \right) 
- \rho\,u_\subN\,\matpartu\,\vec\nu - \nabs\,p 
\nonumber \\ & \quad 
+ \vec g_\subT\,.
\end{align}
We observe that \eqref{eq:3dJOR18finala}, \eqref{eq:JOR18finalb}, 
\eqref{eq:JOR18finalc} is precisely \cite[(3.15)]{JankuhnOR18},
on noting our sign convention for the mean curvature.
\end{remark}

\begin{remark} \label{rem:judge}
We remark that the formulation \eqref{eq:JOR18final} has the advantage of
showing that the PDE for the normal component $u_\subN$ of the velocity is not
of parabolic type. However, other natural properties of the original system
\eqref{eq:esns} are more hidden. For example, the fact that for any surface
$\Gamma$, the system \eqref{eq:esns} with $\vec g= \vec 0$ and $\alpha=0$
is satisfied by $(\Gamma(t), \vec u(\cdot,t), p(\cdot, t)) = (\Gamma, \vec b,
0)$ for any $\vec b \in \bR^3$, is more difficult to read of from 
\eqref{eq:JOR18final}. Compare also with Remark~\ref{rem:discconst}.
\end{remark}

\setcounter{equation}{0} 
\section{Radially symmetric solution} \label{sec:App_truesol}

We construct a solution for a radially expanding sphere of radius $r(t)$,
$\Gamma(t) = r(t)\,\bS^{2} = \{ \vec z \in \bR^3 : |\vec z| = r(t)\}$, 
with $\mat{\mathcal P}\, \vec u = \vec 0$, for the generalized equations
\begin{equation} \label{eq:esnsb}
\rho \, \matpartu\,\vec u 
- \nabs\cdot\mat\sigma =
\alpha\, f_\Gamma\,\vec \nu
- \tfrac\theta2 \rho (\nabs\cdot\vec u)\, \vec u
,\qquad
\nabs\cdot\vec u = b,\qquad \mathcal V = \vec u\cdot\vec\nu,
\end{equation}
where we recall that in our scheme \eqref{eq:qGD} we have either $\theta=0$ or
$\theta=1$.
Then it follows from \eqref{eq:areacons} that
\begin{align*} 
\langle b, 1 \rangle_{\Gamma(t)} &
= \langle1, \nabs\cdot\vec u \rangle_{\Gamma(t)}
= \ddt \mathcal{H}^{2}(\Gamma(t)) 
= \ddt \mathcal{H}^{2}(r(t)\,\bS^{2})
= \mathcal{H}^{2}(\bS^{2})\,\ddt\, [r(t)]^{2} \nonumber \\ &
= 2\,r(t)\,r'(t)\,\mathcal{H}^{2}(\bS^{2})\,,
\end{align*}
and so we can choose the spatially constant function
\begin{equation*} 
b(t) = 2\,\frac{r'(t)}{r(t)}\,.
\end{equation*}
Choosing $\vec\nu = \frac{\vec\id}{|\vec\id|}$ as the outer normal on
$\Gamma(t)$, we have that $\vec u \cdot \vec\nu = \mathcal V = r'(t)$. Hence we
postulate $\vec u(\cdot,t) = r'(t)\,\vec\nu$ as the solution, together with a 
constant pressure $p(t) \in \bR$. On recalling \eqref{eq:diff2} it then holds 
that
\begin{align*}
\nabs\cdot\mat\sigma & 
= 2\,\mu\,\nabs\cdot\mat D_s(\vec u) - \nabs\,p - \varkappa\,p\,\vec\nu
= 2\,\mu\,\nabs\cdot( \vec u \cdot \vec\nu\, \nabs\,\vec\nu) 
- \varkappa\,p\,\vec\nu \nonumber \\ &
= 2\,\mu\,\nabs\cdot( r'(t)\, \nabs\,\vec\nu) - \varkappa\,p\,\vec\nu
= 2\,\mu\,r'(t)\,\Delta_s\,\vec\nu - \varkappa\,p\,\vec\nu
\nonumber \\ &
= -4\,\mu\,\frac{r'(t)}{r^2(t)}\,\vec\nu + 
\frac{2}{r(t)}\,p\,\vec\nu\,,
\end{align*}
where we have noted from Lemmas~12 and 16 in \cite{bgnreview} that on the 
sphere $\Gamma(t)$ we have 
$\Delta_s\,\vec\nu = - |\nabs\,\vec\nu|^2\,\vec\nu= 
-\frac{2}{r^2(t)}\,\vec\nu$, as well as $\varkappa = - \frac{2}{r(t)}$.

We also compute $\rho\,\matpartu\,\vec u = \rho\,r''(t)\,\vec\nu$, so that
overall \eqref{eq:esnsb} will be satisfied if
\begin{equation*} 
\rho\,r''(t)\,\vec\nu 
+4\,\mu\,\frac{r'(t)}{r^2(t)}\,\vec\nu 
- \frac{2}{r(t)}\,p\,\vec\nu = \alpha\,f_\Gamma\,\vec\nu\,.
\end{equation*}
Hence the solution to \eqref{eq:esnsb} is $\Gamma(t) = r(t)\,\bS^{2}$,
$\vec u(\cdot, t) = r'(t)\,\vec\nu$ and
\begin{equation*} 
p(t) = \rho\,\frac{r''(t)\,r(t)}{2} + 2\,\mu\,\frac{r'(t)}{r(t)} 
- \alpha\,\frac{f_\Gamma(t)\,r(t)}{2}\,,
\end{equation*}
where for \eqref{eq:fGamma} we have
\begin{equation*} 
f_\Gamma(t) = \frac{(2)^2}{r^3(t)} - \tfrac12\left(\frac{2}{r(t)}\right)^3
=0\,.
\end{equation*}
Finally, we note that
\[
-\tfrac\theta2\rho(\nabs\cdot\vec u) \vec u 
= -\theta\rho \frac{[r'(t)]^2}{r(t)}\vec\nu
= \tfrac12 \theta\rho [r'(t)]^2 \vec\varkappa,
\]
where we have noted that 
$\vec\varkappa = \varkappa\vec\nu=- \frac{2}{r(t)}\vec\nu$.

\end{appendix}

\def\soft#1{\leavevmode\setbox0=\hbox{h}\dimen7=\ht0\advance \dimen7
  by-1ex\relax\if t#1\relax\rlap{\raise.6\dimen7
  \hbox{\kern.3ex\char'47}}#1\relax\else\if T#1\relax
  \rlap{\raise.5\dimen7\hbox{\kern1.3ex\char'47}}#1\relax \else\if
  d#1\relax\rlap{\raise.5\dimen7\hbox{\kern.9ex \char'47}}#1\relax\else\if
  D#1\relax\rlap{\raise.5\dimen7 \hbox{\kern1.4ex\char'47}}#1\relax\else\if
  l#1\relax \rlap{\raise.5\dimen7\hbox{\kern.4ex\char'47}}#1\relax \else\if
  L#1\relax\rlap{\raise.5\dimen7\hbox{\kern.7ex
  \char'47}}#1\relax\else\message{accent \string\soft \space #1 not
  defined!}#1\relax\fi\fi\fi\fi\fi\fi}

\end{document}